\documentclass{amsart}

\usepackage{a4wide,amssymb}

\usepackage[arrow,matrix,curve]{xy}\SilentMatrices
\def\xyma{\xymatrix@M.7em}

\def\declaration#1#2{
\expandafter\def\csname #1\endcsname##1{\subsection{#2}\label{#1:##1}} }

\declaration{dfn}{Definition}
\declaration{xmp}{Example}
\declaration{thm}{Theorem}
\declaration{prp}{Proposition}
\declaration{crl}{Corollary}
\declaration{lem}{Lemma}
\declaration{rem}{Remark}

\def\ax{\rtimes}

\def\ge{\geqslant}
\def\le{\leqslant}

\def\ph{\varphi}

\def\gth{\mathfrak}
\def\eps{\varepsilon}
\def\then{\Rightarrow}
\def\set#1{\left\{#1\right\}}
\def\setof#1#2{\left\{#1\!\mid\!#2\right\}}
\def\xto#1#2{\xrightarrow[#2]{#1}}
\def\brk#1{\left\langle{#1}\right\rangle}

\def\x{\operatornamewithlimits{\times}}
\def\o{\circ}
\def\os{\operatorname{\oplus}}
\def\ox{\operatornamewithlimits{\otimes}}
\def\op{^\mathrm{op}}

\def\hom{\operatorname{hom}}
\def\Hom{\operatorname{Hom}}

\def\im{\operatorname{Im}}
\def\ker{\operatorname{Ker}}
\def\ext{\operatorname{Ext}}
\def\coker{\operatorname{Coker}}

\def\id{\operatorname{id}}
\def\toto{\rightrightarrows}
\def\into{\rightarrowtail}
\def\onto{\twoheadrightarrow}
\def\d{\partial}
\def\inv{^{-1}}

\def\s{\mathbb S}
\def\t{\mathbb T}
\def\a{\mathbb A}
\def\b{\mathbb B}

\def\e{\mathbb E}
\def\1{{1\kern-.22em\mathrm l}}
\def\z{\mathbf Z}
\def\A{\mathcal A}
\def\co{\mathcal C}
\def\C{\mathbf C}
\def\V{\mathbf V}

\def\Set{\mathbf{Set}}
\def\cat{\mathbf{Cat}}
\def\ab{\mathbf{Ab}}
\def\Th{\mathbf{Theories}}
\def\Lf{\mathbf{Lf}}
\def\m{\textrm-\mathbf{mod}}
\def\aff{\textrm-\mathbf{aff}}
\def\tup#1{\left<#1\right>}
\def\der{\operatorname{Der}}
\def\ider{\operatorname{Ider}}
\def\ad{\operatorname{ad}}
\def\gd{\delta}
\def\mic{\Delta}
\def\mac{\nabla}
\def\sle{/\!\!_{_=}}

\begin{document}

\title{Linear extensions and nilpotence of Maltsev theories}

\author{Mamuka Jibladze}
\author{Teimuraz Pirashvili}
\address{Department of Algebra, Razmadze Mathematical Institute, Tbilisi
380093, Georgia} 

\thanks{Research supported by the grant INTAS-99-00817 and the TMR network
``K-theory and algebraic groups'', ERB FMRX CT-97-0107}

\begin{abstract}
Relationship is clarified between the notions of linear extension of algebraic theories, and central extension, in the sense of commutator calculus, of their models. Varieties of algebras turn out to be nilpotent Maltsev precisely when their theories may be obtained as results of iterated linear extensions by bifunctors from the so called abelian theories. The latter theories are described; they are slightly more general than theories of modules over a ring.
\end{abstract}

\maketitle

\section*{Introduction}

The notion of linear extension of categories introduced in \cite{bw} is
basic in the study of cohomological properties of algebraic theories.
Roughly, linear extensions play the same r\^ole for theories as extensions
with abelian kernel for groups. It seems that many remarkable properties
of theories are preserved under linear extensions. As an example, one can
mention the fact (proved in \cite{pirashvili}) that if a theory has the
property that all of its projective models are free, then the same is true
for any other theory obtained from it by linear extensions.

One of the goals of the present paper is to investigate behaviour of
Maltsev theories (the ones possessing a ternary operation $p$ which is 
Maltsev, i.~e. obeys the identity $p(x,x,y)=p(y,x,x)=y$) under linear
extensions. In particular, it turns out (Proposition \ref{prp:mal}) that
any linear extension of a Maltsev theory is itself Maltsev.

Another aspect of the r\^ole that linear extensions can play in the study of algebraic theories is related to the notion of nilpotence for theories. One might ask whether there is an analog for theories of the fact that nilpotence of an
algebra (say, Lie algebra, or a group, or an associative algebra without unit) is equivalent to the existence of a finite tower of central extensions starting with an abelian algebra and ending with the given algebra. So one might call a theory nilpotent if there is a tower
of linear extensions of theories starting with a theory which is
``abelian'' and ending with the given one. Such definition is inherent in
\cite{jp}, where it is proved that algebraic theories corresponding to the
varieties of nilpotent groups (resp.  algebras) of nilpotence class $n$
fit into towers of linear extensions of length $n$ by certain bifunctors,
starting from abelian, or linear, theories, i.~e. theories of modules over
some ring. See also \cite{pirashvili}.

On the other hand, there is a well understood
generalization of the commutator calculus from groups or Lie algebras to
much more general varieties of universal algebras (the initial idea is
contained in \cite{smith}; a maximally exhaustive treatment is probably
\cite{fm}). In particular there is a notion of abelian (linear) and
nilpotent varieties, generalizing the ones for groups and algebras. This approach yields most satisfactory results for Maltsev varieties. One can ask, what is the relationship between these two approaches.

It will be proved that these approaches are indeed equivalent for Maltsev
theories. That is, a Maltsev variety is nilpotent of class $n$ in the sense
of commutator calculus if and only if the corresponding theory can be
obtained by $n$-fold linear extensions of particular ``untwisted'' type,
starting from a Maltsev theory which is abelian in the sense of commutator
calculus. In a section by the first author, a description of such abelian
Maltsev theories and linear extensions between them is given.

The first author gratefully acknowledges discussions with G.~Janelidze
which helped him to realize the r\^ole of torsors under constant groups
for characterizing centrality of an abelian extension, as in
\ref{prp:malcen} below.

Concerning notation: throughout the paper, it is of set-theoretic style, as
is by now usual in category theory. For example, for a morphism $f:G\to H$
between internal groups in a category with products, an expression like
$f(x+y)$ might mean composite of $f$ with $+(x,y)$, for variable morphisms
$x,y:X\to G$, with $+:G\x G\to G$ being the group operation. Or, for a
congruence on an object $A$, i.~e. a parallel pair $(a_1,a_2):R\toto A$ such
that the resulting map $\tup{a_1(-),a_2(-)}:\hom(X,R)\to\hom(X,A)\x\hom(X,A)$
is an equivalence relation for any $X$, we might use notation $xRy$ or
$x\sim_Ry$, or just $x\sim y$ for morphisms $x,y:X\to A$ such that $(x,y)$
factors through $R$. 

\section{Linear extensions as torsors}

Various kinds of extensions that appear in this paper are based on the notion of a \emph{principal $G$-object} or \emph{$G$-torsor}, or \emph{torsor under} $G$, for an internal group $G$ in a category $\C$ with finite products. It is a $G$-object $T$ which
satisfies two conditions; the first condition is that the morphism (action, projection) : $G\x T\to T\x T$ is an isomorphism. The second
condition says that $T$ has global support and may have different meanings,
depending on exactness properties of $\C$. In this paper this condition
usually means that $T\x T\toto T\to1$ is a coequalizer.

The first condition is often expressed equationally, using a ``subtraction'' map $-:T\x T\to G$, namely the composite of the projection $G\x T\to G$ with the inverse of the above isomorphism. It is easy to see that the condition is equivalent to requiring two identities
\begin{align*}
(g+x)-x&=g,\\
(x-y)+y&=x
\end{align*}
for $g:X\to G$, $x,y:X\to T$. In fact the whole torsor structure can be
expressed by requiring that the morphism $m:T\x T\x T\to T$ given by
$m(x,y,z)=(x-y)+z$ be an associative Maltsev operation, in the following
sense: 

\dfn{maltsev}
A morphism $m:T\x T\x T\to T$ is called a \emph{Maltsev operation} on $T$
if it satisfies
$$
m(x,y,y)=x=m(y,y,x).
$$
It is called \emph{associative} if one has
$$
m(u,v,m(x,y,z))=m(m(u,v,x),y,z)
$$
and \emph{commutative} if $m(x,y,z)=m(z,y,x)$.

\

Note for future reference:

\lem{asmal}
{\sl
Any associative Maltsev operation $m$ satisfies}
$$
m(u,v,m(x,y,z))=m(u,m(y,x,v),z).
$$
\begin{proof}
For readability, denote $m(x,y,z)=x-y+z$. We thus have
$$
x-y+y=x=y-y+x
$$
and
$$
u-v+(x-y+z)=(u-v+x)-y+z.
$$
We can thus save parentheses and denote the latter expression by $u-v+x-y+z$.
One then has 
\begin{align*}
    u-v+x-y+z &= (u-(y-x+v)+(y-x+v))-v+x-y+z\\
  &= u-(y-x+v)+(y-x+v-v+x-y+z)\\
  &= u-(y-x+v)+z.
\end{align*}
\end{proof}

One has the following (well known, in various guises) fact:

\prp{maltor}
{\sl
Let $\C$ be a category with coequalizers of congruences and finite products,
which commute (that is, coequalizer of a product of diagrams is product of
their coequalizers). Then there is a torsor structure on an object $T$ iff 
it has global support and there is an associative Maltsev operation on $T$. In
fact there is a one-to-one correspondence between such structures. Moreover
commutative Maltsev operations correspond to structures of torsors under
abelian groups.}
\begin{proof}
Given a torsor structure, the Maltsev conditions for $m(x,y,z)=(x-y)+z$ are
just the two identities above. As for associativity, it means
$$
(u-v)+((x-y)+z)=(((u-v)+x)-y)+z,
$$
which follows easily from
$$
(g+x)-y=(g+((x-y)+y))-y=((g+(x-y))+y)-y=g+(x-y).
$$
If moreover the group is commutative, one has
$$
(x-y)+z=(x-y)+(z-y)+y=(z-y)+(x-y)+y=(z-y)+x.
$$

Conversely, for an associative Maltsev operation $m$ on $T$, the relation
$$
(x,y)\sim(m(x,y,z),z)
$$
is a congruence on $T\x T$. Indeed, it is reflexive since
$(m(x,y,y),y)=(x,y)$, symmetric since $(m(x,y,z),z)$ $\sim$
$(m(m(x,y,z),z,y),y)$ $=$ $(m(x,y,m(z,z,y)),y)$ $=$ $(x,y)$, and transitive
since $(m(m(x,y,z),z,t),t)$ $=$ $(m(x,y,m(z,z,t)),t)$ $=$ $(m(x,y,t),t)$. Let
$G$ be the coequalizer, and let $-:T\x T\onto G$ be the quotient map. Since
$(m(x,x,y),y)=(y,y)$, one has $(x,x)\sim(y,y)$, i.~e. $x-x=y-y$. In
particular, taking for $x,y$ the projections $T\x T\to T$ we get a map from
their coequalizer to $G$, i.~e. a global element $0:1\to G$, by the global
support condition on $T$. Addition on $G$ is defined by
$(x-y)+(z-t)=m(x,y,z)-t$ which is legitimate since cartesian product of two
coequalizers is a coequalizer in our category. Then additive inverse of $x-y$
is $y-x$, and the action of $G$ on $T$ is given by $(x-y)+z=m(x,y,z)$. It is
straightforward to verify all the remaining identities. 
\end{proof}

So torsors in such categories determine (at least one of) their own groups.
In view of this, objects equipped with an associative Maltsev operation will
be also referred to as torsors. In the literature they are also known as
\emph{herds}. 

\dfn{abcen}
A morphism $p:E\to B$ in a category $\C$ with products is called an
\emph{abelian extension}, or simply abelian, if it admits a structure of a
torsor in $\C/B$, for some internal abelian group in $\C/B$. If furthermore
the group has the form $B^*(A)=(B\x A\to B)$, for some internal abelian group
$A$ in $\C$, then the morphism is called \emph{central extension}. 

\

So in view of \ref{prp:maltor}, a morphism $p:E\to B$ in a sufficiently nice
category is an abelian extension iff it is coequalizer of its own kernel pair
$E\x_BE\toto E$, and there is an associative Maltsev operation $E\x_BE\x_BE\to
E$ over $B$. As for central extensions, one has another omnipresent fact: 

\prp{malcen}
{\sl
A morphism $p:E\to B$ in $\C/B$ is a torsor under a constant group $B^*(G)$,
for some group $G$ in $\C$, if and only if the corresponding Maltsev operation
$m:E\x_BE\x_BE\to E$ extends to an associative Maltsev operation $E\x_BE\x
E\to E$, with $pm(x,y,z)=pz$ for any $(x,y,z):X\to E\x_BE\x E$.}
\begin{proof}
If the group is $B^*(G)$, then the action can be written as $+:G\x E\cong(B\x
G)\x_BE\to E$, and the subtraction is given by $(p,-):E\x_BE\to B\x G$ for
some map $-:E\x_BE\to G$. Thus $(x-y)+z$ is defined for $px=py$ and any $z$.
Moreover using $p(g+x)=px$, all the identities are proved in exactly the same
way as in \ref{prp:maltor}. 

Conversely if $m$ as above is given, we construct the group $G$ as quotient
of $E\x_BE$ by $(x,y)\sim(m(x,y,z),z)$ again, with $(x,y):X\to E\x_BE$ and
\emph{any} $z:X\to E$. That is, we coequalize the maps $(x,y,z)\mapsto(x,y)$
and $(x,y,z)\mapsto(m(x,y,z),z)$ from $E\x_BE\x E$ to $E\x_BE$. This then
gives the maps $+:G\x E\to E$ and $-:E\x_BE\to G$ just as in \ref{prp:maltor},
satisfying the required identities. 
\end{proof}

Terminology in \ref{dfn:abcen} above is motivated by one important case, when
$\C$ is the slice $\V/B$ of some variety $\V$ of universal algebras over one
of its objects $B$. Then the object $p:E\to B$ of $\V/B$ having global support
simply means that $p$ is surjective. It is well known that this notion of
torsor gives various kinds of extensions of universal algebras. For example,
when $\V$ is the variety of groups, then for any internal group $G$ in $\V/B$
there is a $B$-module $M$ such that $G$ is isomorphic to the projection
$B\ax M\to B$ of the semidirect product of $B$ with $M$, the group structure
given by homomorphisms ($+:(B\ax M)\x_B(B\ax M)\to B\ax M$,
$-:B\ax M\to B\ax M$, $0:B\to B\ax M$) with
$((b,x_1),(b,x_2))\mapsto(b,x_1+x_2)$, $(b,x)\mapsto(b,-x)$, $b\mapsto(b,0)$
respectively. Furthermore a $G$-torsor is the same as a short exact sequence
$$
\xyma{M\ar@{>->}[r]^i&E\ar@{->>}[r]^p&B}
$$
with $i(p(e)x)=e+i(x)-e$, the action $(B\ax M)\x_BE\to E$ being given by
$((b,x),e)\mapsto i(x)+e$. Similarly when $\V$ is the variety of Lie rings,
an internal group in $\V/B$ amounts to a $B$-module $M$, and a torsor under
this group to an extension $0\to M\to E\to B\to0$; when $\V$ is the variety
of associative algebras with unit, one gets bimodules and singular extensions,
etc. Note that in all these cases torsors under \emph{constant} groups, i.~e.
internal groups in $\V/B$ represented by projections $B\x G\to B$ for an
internal group $G$ in $\V$, correspond to central extensions of $B$.

Another situation where torsors are important for us arises from a small
category $\b$, with $\C$ the full subcategory $\cat\sle\b$ of the slice
$\cat/\b$ of categories over $\b$, consisting of those functors
$p:\e\to\b$ which are identity on objects. In this case the global support
condition is that $p$ is full, i.~e. surjective on morphisms. The resulting
notion turns then out to be equivalent to the notion of \emph{linear
extension} of categories, which we now recall.

For a small category $\b$, let $\b^\#$ denote the category called
\emph{twisted arrow category} of $\b$ in \cite{maclane}, and the
\emph{category of factorizations} of $\b$ in \cite{bw}. Objects of
$\b^\#$ are morphisms of $\b$, whereas $\hom_{\b^\#}(b,b')$ consists of
pairs $(b_1,b_2)$ with $b_1bb_2=b'$. A \emph{natural system} on $\b$
with values in a category $\C$ is a functor $D:\b^\#\to\C$. It is thus a
collection of $\C$-objects $(D_b)_{b:X\to Y}$ of $\C$, indexed by morphisms of
$\b$, together with $\C$-morphisms $b_1(\ ):D_b\to D_{b_1b}$ and
$(\ )b_2:D_b\to D_{bb_2}$, for all composable morphisms $b_1$, $b$, $b_2$ in
$\b$, such that certain evident diagrams commute. In other words, one must
have 
\begin{align*}
(b_1b_2)x_3&=b_1(b_2x_3),\\
(b_1x_2)b_3&=b_1(x_2b_3),\\
(x_1b_2)b_3&=x_1(b_2b_3)
\end{align*}
for any composable $\xto{b_3}{}\xto{b_2}{}\xto{b_1}{}$ and any
$x_i:X\to D_{b_i}$. 

We will use the following notion from \cite{bw}: for a natural system $D$ on
a category $\b$ with values in abelian groups, a \emph{linear extension}
of $\b$ by $D$ is an object of $\cat\sle\b$, i.~e. a functor
$P:\e\to\b$ that is identity on objects, together with transitive and
effective actions
$D_b\x P\inv(b)\to P\inv(b)$, $(x,e)\mapsto x+e$, for all $b:X\to Y$ in $\b$,
such that for any composable morphisms $e_1$, $e_2$ in $\e$ and any
$x_i\in D_{P(e_i)}$, $i=1,2$, one has
$$
(x_1+e_1)(x_2+e_2)=(x_1P(e_2)+P(e_1)x_2)+e_1e_2.
$$

An example of a linear extension by a natural system $D$ is given by the
\emph{trivial} linear extension $\b\ax D$ with
$\hom_{\b\ax D}(X,Y)=\coprod_{b:X\to Y}D_b$, composition
$x_1x_2=x_1b_2+b_1x_2$ for $x_1\in D_{b_1}$, $x_2\in D_{b_2}$, identities
$0\in D_{1_X}$, $X\in\b$, and the actions $D_b\x D_b\to D_b$ given by the
group law in $D_b$. 

For natural systems $D$ of abelian groups, there are cohomology groups
$H^*(\b;D)$ of $\b$ with coefficients in $D$, having the usual properties,
such that $H^2(\b;D)$ classifies linear extensions of $\b$ by $D$.
See \cite{bw}.

\prp{lin=tor}
{\sl
For any small category $\b$, assigning to a natural system $D$ on $\b$ the
trivial extension $\b\ax D\to\b$ determines an equivalence between the
category of natural systems of abelian groups on $\b$ and the category of
internal abelian groups in $\cat\sle\b$. Moreover linear extensions of $\b$
are the same as torsors in $\cat\sle\b$, i.~e. those objects which, as
morphisms in $\cat$, are abelian in the sense of \ref{dfn:abcen}; more
precisely, for any natural system $D$ on $\b$ linear extensions of $\b$ by $D$
are in one-to-one correspondence with $(\b\ax D\to\b)$-torsors in
$\cat\sle\b$.}
\begin{proof}
The group structure on $\b\ax D$ is given as follows: the zero is the functor
$\b\to\b\ax D$ which sends a morphism $b:X\to Y$ to $0\in D_b$. The addition
functor $(\b\ax D)\x_\b(\b\ax D)\to\b\ax D$ is given by addition in the groups
$D_b$ and similarly for inverses. Conversely, any abelian group $P:\a\to\b$ in
$\cat\sle\b$ determines a natural system with $D_b=P\inv(b)$. These
correspondences are evidently functorial and can be easily checked to
define mutually inverse equivalences.

Similarly, given any linear extension $\e\to\b$, of $\b$ by a natural
system $D$, its transitive and effective actions combine into a functor
$(\b\ax D)\x_\b\e\to\e$ which can be checked to form a $\b\ax D$-torsor.
And conversely, any torsor furnishes the required action for a linear
extension.
\end{proof}

In particular, linear extensions can be defined in terms of the subtraction
map. One sees easily that the corresponding identities are
\begin{equation}
e_1e_2-e_1e_2'=P(e_1)(e_2-e_2'),\ \ e_1e_2-e_1'e_2=(e_1-e_1')P(e_2),
\tag{$\bullet$}
\end{equation}
for $P(e_i)=P(e_i')$, $i=1,2$. In view of \ref{prp:maltor} and
\ref{prp:lin=tor}, linear extensions can be also defined in terms of
commutative associative Maltsev operations, without mentioning any natural
system. Namely, linear extension structures on an object $\e\to\b$ of
$\cat\sle\b$ with global support are in one-to-one correspondence with
functors $\e\x_\b\e\x_\b\e\to\e$ over $\b$ which are commutative associative
Maltsev operations. 

There is another context in which natural systems arise as internal abelian
groups. 

\prp{bifuslice}
{\sl
For any $\b$, there is an equivalence of categories
$$
\ab(\cat\sle\b)\cong\ab(\Set^{\b\op\x\b}/\hom_\b),
$$
i.~e. the category of natural systems of abelian groups on $\b$ is
equivalent to the category of internal abelian groups in the slice
over $\hom_\b$ of the category of set-valued bifunctors on $\b$.
Under this equivalence, the inclusion
$$
\hom_\b^*:\ab(\Set^{\b\op\x\b})\to\ab(\Set^{\b\op\x\b}/\hom_\b)
$$
of constant internal groups, carrying $D$ to the projection
$\hom_\b\x D\to\hom_\b$, becomes identified with the inclusion into
natural systems, carrying a bifunctor $D:\b\op\x\b\to\ab$ to the
natural system with $D_{b:X\to Y}=D(X,Y)$.
} \qed

This (as well as \ref{prp:lin=tor}, in fact) is a consequence of general facts
from \cite{bjt} (see 1.5 and 4.11 there). 
 
We will call the particular natural systems arising, as above, from bifunctors,
and linear extensions by them \emph{untwisted}. 

Thus natural systems on $\b$ can be identified with (trivial) abelian
extensions of $\hom_\b$ in $\Set^{\b\op\x\b}$, in the sense of
\ref{dfn:abcen}, and moreover untwisted natural systems correspond to
trivial central extensions in that category. A natural question then
arises --- what should be analog of \ref{prp:malcen} in this context, that
is, which torsors in $\cat\sle\b$ correspond to untwisted linear
extensions under the equivalence of \ref{prp:lin=tor}. The answer is given
by the following

\prp{simtor} {\sl Let $P:\e\to\b$ be a full functor bijective on objects,
with a torsor structure in $\cat\sle\b$ given by the functor
$m:\e\x_\b\e\x_\b\e\to\e$ over $\e$. Then, the linear extension
corresponding to it by \ref{prp:lin=tor} is untwisted if and only if $m$
can be extended to a collection of commutative associative Maltsev
operations $$
m_{X,Y}:\hom_\e(X,Y)\x_{\hom_\b(X,Y)}\hom_\e(X,Y)\x\hom_\e(X,Y)\to\hom_\e(X,Y),
$$ such that $$ P(m_{X,Y}(f_1,f_2,f))=P(f)  $$ and $$
gm_{X,Y}(f_1,f_2,f)=m_{X,Z}(gf_1,gf_2,gf),\ \
m_{X,Y}(f_1,f_2,f)h=m_{T,Y}(f_1h,f_2h,fh)  $$ for any
$f_1,f_2,f\in\hom_\e(X,Y)$ with $P(f_1)=P(f_2)$ and any}
$g\in\hom_\e(Y,Z)$, $h\in\hom_\e(T,X)$.

\begin{proof}
The ``only if'' part follows since as soon as one has a linear
extension by a bifunctor $D$, all the groups $D_b$, $D_{b'}$ are naturally
identified for any $b,b'\in\hom_\b(X,Y)$. Hence one can define
$m_{X,Y}(f_1,f_2,f)=(f_1-f_2)+f$ by identifying $f_1-f_2\in D_{P(f_i)}$
with the corresponding element of $D_{P(f)}$. The above identities then
follow easily from the corresponding identities for linear extensions. 

For the ``if'' part, let $D$ be the natural system of abelian groups
corresponding to $P$ according to \ref{prp:maltor} and \ref{prp:lin=tor}. Thus for $b\in\hom_\b(X,Y)$, elements of $D_b$ have the form $f_1-f_2$, with
$P(f_1)=P(f_2)=b$, and $f_1-f_2$ determines the same element as
$m(f_1,f_2,f)-f$, for any other $f$ with $P(f)=b$. Then using the $m_{X,Y}$
above we can define a collection of isomorphisms $\ph_{b,P(f)}:D_b\to
D_{P(f)}$ for $f\in\hom_\e(X,Y)$, $b\in\hom_\b(X,Y)$, via 
$$
\ph_{b,P(f)}(f_1-f_2)=m_{X,Y}(f_1,f_2,f)-f.
$$
This is correctly defined since
$m_{X,Y}(m(f_1,f_2,f_3),f_3,f)=m_{X,Y}(f_1,f_2,f)$ by associativity and
Maltsev identity. And it does not really depend on $f$. Indeed for any other
$f'$ with $P(f')=P(f)$ one has, by the same identities,
$m(m_{X,Y}(f_1,f_2,f),f,f')$ $=$ $m_{X,Y}(f_1,f_2,f')$. But this is the same
as $m_{X,Y}(f_1,f_2,f)-f=m_{X,Y}(f_1,f_2,f')-f'$ -- recall (from the proof of
\ref{prp:maltor} that $m(g,f,f')=g'$ is equivalent to $g-f=g'-f'$ as soon as
$P(f)=P(g)=P(f')=P(g')$.

Furthermore,
$$
\ph_{b,b}=\textrm{identity}_{D_b},
$$
as $m_{X,Y}(f_1,f_2,f)-f=m(f_1,f_2,f)-f=f_1-f_2$ for $P(f_1)=P(f_2)=P(f)$, and
$$
\ph_{b',b''}\ph_{b,b'}=\ph_{b,b''},
$$
as $m_{X,Y}(m_{X,Y}(f_1,f_2,f'),f',f'')-f''=m_{X,Y}(f_1,f_2,f'')-f''$ for
$P(f_1)=P(f_2)=b$, $P(f')=b'$, $P(f'')=b''$.

Finally for any $a:Y\to Z$, $b,b':X\to Y$ in $\b$, and any $f_1$, $f_2$ with
$P(f_1)=P(f_2)=b$ one has, using the equations ($\bullet$) above:
\begin{multline*}
a\ph_{b,b'}(f_1-f_2)=a(m_{X,Y}(f_1,f_2,f)-f)=gm_{X,Z}(f_1,f_2,f)-gf\\
=m_{X,Z}(gf_1,gf_2,gf)-gf=\ph_{ab,ab'}(gf_1-gf_2)=\ph_{ab,ab'}(a(f_1-f_2))
\end{multline*}
for any $g$ with $P(g)=a$, and similarly, for any $c:T\to X$ in $\b$,
$$
(\ph_{b,b'}(f_1-f_2))c=\ph_{bc,b'c}((f_1-f_2)c).
$$
We thus have constructed an isomorphism of the natural system $D$ with the one
obtained from the bifunctor $\bar D$, where
$$
\bar D(X,Y)=\left(\bigoplus_{b\in\hom_\b(X,Y)}D_b\right)/\sim,
$$
with $\sim$ being the equivalence relation identifying any $x\in D_b$ with
$\ph_{b,b'}(x)\in D_{b'}$, for all $b,b'\in\hom_\b(X,Y)$.
\end{proof}

\section{Theories}

\subsection*{Recollections on algebraic theories}
Everywhere in the sequel, $\Set$ will denote the category of sets. The
opposite of its full subcategory with finite sets $\set{1,...,n}$ ($n\ge0$) as
objects, will be denoted $\s$. Its objects will be redenoted by $X^0=1$,
$X^1=X$ (the \emph{generator}), $X^2$, $X^3$, ... , and the morphisms from
$\s(X^n,X)$ by $x_1,...,x_n$.

A \emph{finitary algebraic theory}, or simply theory, is a small category $\t$
equipped with a functor $\s\to\t$ which is identity on objects and preserves
finite products. This functor will be usually suppressed from the notations,
and objects and morphisms of $\s$ will be identified with their images under
it --- an usual abuse of notation with algebras.

A \emph{model} of a theory $\t$ in a category $\C$ is a finite product
preserving functor from $\t$ to $\C$. These functors and their natural
transformations form the category of models $\t\m(\C)$. For $\C=\Set$, this
will be abbreviated to just $\t\m$. Since representable functors preserve
any available limits, there is a full embedding $I_\t:\t\op\to\t\m$.

A model $M$ of $\t$ is in fact nothing but an object $M(X)$ with operations
$u_M:M(X)^n\to M(X)$ for each element of $u\in\hom_\t(X^n,X)$, satisfying
identities prescribed by category structure of $\t$. By this reason, elements
of $\t(n)=\hom_\t(X^n,X)$ will be called $n$\emph{-ary operations of} $\t$. 
Thus for any theory $\t$, the category $\t\m$ is a variety of universal
algebras. In particular, $\t\m$ is an exact category, regular epis are exactly
surjective maps, etc. Conversely, for any variety $\V$, the opposite of the
category of the algebras freely generated by the sets $\set{1,...,n}$,
$n\ge0$, is an algebraic theory, whose category of models is equivalent to
$\V$. 

A morphism of theories $\t'\to\t$ is a model of $\t'$ in $\t$ which respects
the structure functors from $\s$. So by definition, $\s$ is the initial object
of the category $\Th$ of finitary algebraic theories. For every morphism of
theories $F:\t'\to\t$, the induced forgetful functor ``compose with $F$'',
$U_F:\t\m\to\t'\m$, has a left adjoint, which we again denote by $F$,
with the adjunction unit $\eta:$ Identity $\to (-)_F:=U_FF$ and counit
$\eps:FU_F\to$ Identity. In particular, for $\t'=\s$, the corresponding
adjoint pair will be denoted $\t(-):\Set\rightleftarrows\t\m:U_\t$; models in
the image of $\t(-)$ are \emph{free}. The notation is justified by the fact
that these left adjoints are compatible with the above embeddings in the sense
that for any $F:\t'\to\t$, the square 
$$
\xyma{
{\t'}\op      \ar[r]^F
\ar[d]^{I_{\t'}}            &\t\op
                             \ar[d]^{I_\t}\\
\t'\m        \ar[r]^F       &\t\m
}
$$
commutes, so taking $\t'=\s$, the Yoneda embedding $I_\t$ identifies $\t\op$
with the full subcategory of $\t\m$ consisting of free models generated by
objects of $\s$, i.~e. by finite cardinals. Moreover note that the functor
$U_\t$ is representable by the free model on one generator,
$U_\t(M)\cong\hom(\t(1),M)$ for any $\t$-model $M$. So in the sequel, we will
interchangeably use notation $\t(n)$ for $\hom_\t(X^n,X)$, for $I_\t(X^n)$
and for $\t(\set{1,...,n})$.

We note for future reference the following (doubtlessly well known) fact:

\prp{subvar}
{\sl
A morphism of theories $P:\t'\to\t$ is a full functor if and only if each
$\eta_M:M\to M_P$ is surjective. In this case, $U_P$ is full and faithful,
and the corresponding full replete image of $\t\m$ under $U_P$ is a subvariety
of $\t'\m$, i.~e. a full subcategory closed under subobjects, products and
homomorphic images. Moreover, for any surjection $q:M\onto N$ in $\t'\m$, the
square
$$
\xyma
{
M\ar@{->>}[r]^q\ar@{->>}[d]^{\eta_M}&N\ar@{->>}[d]^{\eta_N}\\
M_P\ar@{->>}[r]^{q_P}&N_P
}
$$
is pushout.}
\begin{proof}
Consider the maps $\hom_{\t'}(X^n,X)\to\hom_\t(X^n,X)$, $n\ge0$, induced by
$P$. These are surjective iff $P$ is full, and can be identified with 
$\t'$-homomorphisms $P_n:I_{\t'}(n)\to I_\t(n)$, namely with
$U_{\t'}\eta_{\t'(n)}:U_{\t'}\t'(n)\to U_{\t'}(\t'(n)_P)=U_\t\t(n)$. 
Thus $P$ is full iff $\eta_{\t'(n)}$ is surjective for each $n$.

Next consider a free $\t'$-model $\t'(S)$, for some set $S$. Then for any
$(s_1,...,s_n):\set{1,...,n}\to S$ there is a commutative square
$$
\xyma{
\t'(n)\ar[rr]^{\t'(s_1,...,s_n)}\ar[d]^{\eta_{\t'(n)}}&&\t'(S)\ar[d]^{\eta_{\t'(S)}}\\
 \t(n)\ar[rr]^{\t(s_1,...,s_n)}&&\t(S),
}
$$
so the homomorphism $\eta_{\t'(S)}$ is given by colimit of a filtered diagram
of homomorphisms of the form $\eta_{\t'(n)}$. As each of these homomorphisms
is surjective, $\eta_{\t'(S)}$ is surjective too iff $P$ is full. Finally for
any model $M$ one chooses a surjective homomorphism $q:\t'(S)\onto M$. Both
$P$ and $U_P$ preserve surjections: the former -- since surjections are
precisely coequalizers, and since $P$, being a left adjoint, preserves them;
and the latter -- since it commutes with the forgetful functors to sets. Thus
$\eta_Mq=q_P\eta_{\t'(S)}$ is surjective, so all $\eta_M$'s are surjective iff
$P$ is full. 

Now by adjunction, $\eta U_P$ is a split mono, with left inverse $U_P\eps$.
In our case $\eta U_P$ is also componentwise surjective, so these natural
transformations are mutually inverse isomorphisms. As $U_P$ obviously reflects
isomorphisms, it follows that $\eps$ is an isomorphism, i.~e. $U_P$ is full
and faithful. 

Next, consider a mono $i:M\into U_P(N)$. Its composite with the isomorphism
$\eta_{U_P(N)}$ is mono again. On the other hand
$\eta_{U_P(N)}i=U_PP(i)\eta_M$, so $\eta_M$ is also mono. As it is surjective,
it is thus an isomorphism, i.~e. $M$ also belongs to the replete image of
$U_P$. 

Now to show that the image of $U_P$ is also closed under quotients, suppose
given a surjective homomorphism $q:U_P(N)\onto M$ in $\t'\m$. For $M$ to
belong to the image of $U_P$, for any $u,u'$ in $\t'(n)$ with
$P_n(u)=P_n(u')$, the maps $M(u),M(u'):M(X^n)\to M(X)$ must be equal. In fact
since the $P_n$ are surjective, this condition is also sufficient. But
$M(u)q(X^n)=q(X)N(u)=q(X)N(u')=M(u')q(X^n)$, and since $q$ is surjective,
$M(u)=M(u')$. 

Finally, let us prove the pushout property of the square above. Indeed, given
homomorphisms $h:N\to N'$, $r:M_P\to N'$ with $hq=r\eta_M$, one has
$\im(h)=\im(hq)=\im(r\eta_M)=\im(q)$. The image of $U_P$ is closed under
quotients, as we just proved, so any quotient of $M_P=U_PP(M)$ is (isomorphic
to) $U_P$ of something. In particular so is $\im(h)$, and by adjunction $h$
factors through $\eta_M$, via some $h'$. Then $h'q_P=r$ since their
composites with the epi $\eta_M$ are easily seen to be equal. 
\end{proof}

\subsection*{Remark}
In fact, conditions of the above proposition are interrelated: it follows, for
example, from 3.1 in \cite{jk1}, that the image of $U_P$ is closed under
subobjects iff $\eta$ is surjective, and under quotients iff the indicated
squares are pushouts. 

\ 

\subsection*{Extensions of theories}
Since for a theory $\t$ the category $\Th/\t$ is a subcategory of
$\cat\sle\t$, closed under finite products, internal groups and torsors in
$\Th/\t$ are particular groups and torsors in $\cat\sle\t$, hence by
\ref{prp:lin=tor} they can be considered as particular natural systems and
linear extensions over $\t$. It is easy to identify the property of
natural systems which distinguishes these particular ones (see \cite{jp}):

\dfn{cartesian}
A natural system $D$ on a category with finite products $\t$ is said to be
\emph{cartesian} if for any product diagram
$p_i:X_1\x...\x X_n\to X_i$, $i=1,...,n$ and any $f:X\to X_1\x...\x X_n$,
the maps $p_i(\ ):D_f\to D_{p_if}$, $i=1,...,n$ also form a product
diagram.

\

The equivalence of \ref{prp:bifuslice} restricted to cartesian natural
systems yields

\prp{cartecoef}
{\sl
The category of cartesian natural systems of sets on a theory $\t$ is
equivalent to the category
$$
\t\m(\Set^{\t\op})/I_\t,
$$
with untwisted cartesian natural systems corresponding to objects in the image
of}
$$
I_\t^*:\t\m(\Set^{\t\op})\to\t\m(\Set^{\t\op})/I_\t.
$$
\begin{proof}
Indeed, looking at the equivalence in \ref{prp:bifuslice} one sees that
the category of cartesian natural systems of sets on a small category with
finite products $\t$ is equivalent to the full subcategory of
$\Set^{\t\op\x\t}/\hom_\t$ on those natural transformations $p:B\to\hom_\t$
for which the natural system given by $b\mapsto p\inv(b)$ is cartesian. But
it is straightforward to check that this happens iff $B$ preserves finite
products in the second variable. Thus when $\t$ is a theory this means that
for any fixed object $X^n$, the functor $B(X^n,-)$ is a model of $\t$. So
cartesian natural systems correspond to the full subcategory
$(\t\m)^{\t\op}\cong\t\m(\Set^{\t\op})$ of
$(\Set\vphantom)^\t)^{\t\op}\cong\Set\vphantom)^{\t\op\x\t}\cong(\Set\vphantom)^{\t\op})^\t$.
\end{proof}

\crl{ab}
{\sl
Every linear extension $P:\t\to\t_R$ of the theory $\t_R$ of (left) modules
over a ring $R$ is untwisted.}
\begin{proof}
The category $\t_R\m=R\m$ is abelian, hence so is $\t_R\m(\Set^{\t_R\op})$ $=$
$(\t_R\m)^{\t_R\op}$; but as it is well known, for any additive category $\A$,
and any of its objects $A$, the functors
$$
\xyma
{
\A&\ab(\A)\ar[l]_-{\textrm{forget}}\ar[r]^-{A^*}&\ab(\A/A)
}
$$
are equivalences of categories. In our case this gives that every object of
$\ab((\t_R\m)^{\t_R\op}/I_{\t_R})$ is isomorphic to a projection
$I_{\t_R}\x T\to I_{\t_R}$ for some $T:\t_R\op\to\t_R\m$. Translating this
fact along the equivalence of \ref{prp:bifuslice} one obtains that any
cartesian abelian natural system on $\t_R$ is isomorphic to one of the form
$D_{u:X^m\to X^n}\equiv\Hom_R(R^n,T(X^m))$, for some functor
$T:\t_R\op\to R\m$. Evidently this is an untwisted natural system.
\end{proof}

One has (cf. \cite{jp}, (6.1))

\prp{cartesian}
{\sl
A natural system of abelian groups $D$ on a category with finite products
$\t$ is cartesian iff for any linear extension $P:\t'\to\t$ of $\t$ by $D$,
the category $\t'$ also has finite products, and $P$ preserves them.
} \qed

In particular, linear extensions of an algebraic theory $\t$ by a cartesian
natural system $D$ are morphisms of theories, and equivalence classes of such
extensions form an abelian group isomorphic to $H^2(\t;D)$.

There are lots of examples of linear extensions of theories in \cite{jp}. Let
us mention those which we will encounter in this paper.

\subsection{Examples}\label{exm}
\subsubsection{}\label{mon}
Consider the functor from theories to monoids given by
$\t\mapsto\hom_\t(X,X)$. This functor has a full and faithful right adjoint
assigning to a monoid $M$, the theory $\t_M$ of $M$-sets. Thus the category of
monoids can be identified with a full subcategory of $\Th$ closed under limits
there. In particular, groups, torsors, herds, natural systems, linear
extensions, etc. of monoids (considered as categories with one object) can be
identified with those of the corresponding theories. In other words, a
morphism of theories $P:\t_N\to\t_M$ induced by a homomorphism of monoids
$p:N\to M$ is a linear extension iff $p$, considered as a functor between
categories with one object, is a linear extension -- i.~e. $p$ is an abelian
extension in the category of monoids. The corresponding natural system on $M$
consists of abelian groups $D_x$, for $x$ in $M$, and actions
$x(-):D_y\to D_{xy}$, $(-)y:D_x\to D_{xy}$. It can be also considered as an
``$M$-graded $M$-$M$-bimodule''. The corresponding extensions of theories are
untwisted iff all the $D_x$ are equal.
\subsubsection{}\label{ring}
Any homomorphism of rings $p:S\to R$ gives rise to a morphism $P:\t_S\to\t_R$
from the theory of (left) $S$-modules to that of $R$-modules. This morphism is
a linear extension iff $p$ is a \emph{singular extension}, i.~e. $\ker(p)=B$
is a square zero ideal in $S$. In \cite{jp}, an isomorphism is obtained
$$
H^2(\t_R;D_B)\cong H^2(R;B)
$$
from the group of (untwisted) linear extensions of $\t_R$ by the bifunctor
given by
$$
D_B(X^n,X^k)=\Hom_{R\m}(\t_R(k),B\ox_R\t_R(n))\cong(B\vphantom)^{\os n})^k,
$$
to the second MacLane cohomology group of $R$ with coefficients in $B$.
\subsubsection{}\label{nil}
It is proved in \cite{jp} that for each $n$ there is a linear extension
from the theory of ($n+1$)-nilpotent groups to that of $n$-nilpotent ones;
similarly for groups replaced by Lie rings, associative rings without unit,
or associative commutative rings without unit.
\subsubsection{}\label{mod}
For a left module $M$ over a ring $R$, let $M/(R\m)$ be the coslice category
of modules under $M$, with objects of the form $M\to N$ and obvious
commutative triangles as morphisms. Let $P:M/(R\m)\to R\m$ be the functor
sending $f:M\to N$ to $\coker(f)$. It has a right adjoint $U_P$ given by
$U_P(N)=0:M\to N$. It is then easy to see that this adjoint pair is induced by
a morphism of theories $P:\t_{R;M}\to\t_R$, where $\t_{R;M}$ is the opposite
of the full subcategory of $M/(R\m)$ on objects of the form $(1,0):M\to M\os
R^n$ for $n\ge0$. In particular, $M/(R\m)$ is equivalent to $\t_{R;M}\m$.

Now this $P$ in fact presents $\t_{R;M}$ as a trivial linear extension of
$\t_R$, by the bifunctor $H_M$ given by composition
$$
\t_R\op\x\t_R\xto{\textrm{projection}}{}\t_R\xto{I_{\t_R}\op}{}(R\m)\op\xto{\Hom_R(-,M)}{}\ab,
$$
that is,
$$
H_M(X^n,X^k)=\Hom_R(R^k,M)\cong M^k.
$$
Indeed, the trivial extension $P:\t_R\ax H_M\to\t_R$ can be easily calculated;
one has
$$
\hom_{\t_R\ax H_M}(X^n,X^k)=\Hom_R(R^k,M\os R^n).
$$
One can represent the latter group also as
$$
\hom_{M/(R\m)}(M\xto{(1,0)}{}M\os R^k,M\xto{(1,0)}{}M\os R^n),
$$
which is precisely $\hom_{\t_{R;M}}(X^n,X^k)$.

\subsection*{The Maltsev case}
A \emph{Maltsev theory} is a theory $\t$ for which there is a Maltsev
operation on the generating object $X$. The corresponding variety, i.~e. the
category of models, will be also called Maltsev in this case. It is a
classical result of Maltsev that such varieties are precisely those in which
join of congruences coincides with their composition.

\prp{mal}
{\sl
Let a morphism of theories $P:\t'\to\t$ be a linear extension of $\t$
by a natural system $D$. If $\t$ is a Maltsev theory, then $\t'$ is also
Maltsev.}
\begin{proof}
We will prove that for any $m:X^3\to X$ in $\t'$ such that $P(m)$ is
Maltsev, there is a Maltsev $m':X^3\to X$ with $M(m')=P(m)$.

Let $x_1,x_2,...:X^n\to X$ be the projections. Now
$P(m(x_1,x_2,x_2))$ $=$ $P(x_1)$ and $P(m(x_1,x_1,x_2))$ $=$ $P(x_2)$, so
since $P$ is a $\t\ax D$-torsor, the elements $x_1-m(x_1,x_2,x_2)\in D_{x_1}$
and $x_2-m(x_1,x_1,x_2)\in D_{x_2}$ are defined. Denoting by $x:X\to X$ the
identity, let 
$$
m'=(x_2-m(x_1,x_1,x_2))(x_1,P(m))+(m(x,x,x)-x)P(m)+
(x_1-m(x_1,x_2,x_2))(P(m),x_3)+m.
$$
One then has
\begin{multline*}
m'(x_1,x_1,x_2)\\
=(x_2-m(x_1,x_1,x_2))(x_1,x_2)+(m(x,x,x)-x)x_2+(x_1-m(x_1,x_2,x_2))(x_2,x_2)+m(x_1,x_1,x_2)\\
=(x_2-m(x_1,x_1,x_2))+(m(x_2,x_2,x_2)-x_2)+(x_2-m(x_2,x_2,x_2))+m(x_1,x_1,x_2)=x_2
\end{multline*}
and
\begin{multline*}
m'(x_1,x_2,x_2)\\
=(x_2-m(x_1,x_1,x_2))(x_1,x_1)+(m(x,x,x)-x)x_1+(x_1-m(x_1,x_2,x_2))(x_1,x_2)+m(x_1,x_2,x_2)\\
=(x_1-m(x_1,x_1,x_1))+(m(x_1,x_1,x_1)-x_1)+(x_1-m(x_1,x_2,x_2))+m(x_1,x_2,x_2)=x_1.
\end{multline*}
\end{proof}

\section{Commutators and nilpotence}
Commutator calculus has been extended to general varieties of algebraic
systems by several people -- first by J.~D.~H.~Smith \cite{smith} for Maltsev varieties, then extended to more general cases by Gumm, Hagemann and Herrmann, McKenzie and others (see \cite{fm} for precise information).

Taking the point of view of category theory enables one to make more
apparent the invariant properties of the commutator calculus, i.~e.
properties which do not depend on a particular choice of basic operations
for the algebras of the variety. A good example of such approach is
\cite{pedicchio}. Also in \cite{jk1} a notion of central extension is
derived from abstract categorical version of Galois theory, and it is
shown in \cite{jk2} that central extensions in the sense of commutator
calculus can be described in this way too.

For our paper, categorical reformulation of the commutator calculus given in \cite{pedicchio} is most suitable. Let us recall it briefly.

In a category with kernel pairs and coequalizers one may generate from any pair of morphisms $X\toto M$ a congruence on $M$, as the kernel pair of the coequalizer of this pair. In particular, given two congruences $p',p'':R\toto M$ and $q',q'':S\toto M$, one denotes by
$r',r'':\Delta_{R,S}\to R$ the congruence on $R$ defined to be the kernel pair of the coequalizer of the morphisms (diagonal)$q'$,(diagonal)$q'':S\toto R$. Also let $R'\toto R$ be the kernel pair of $p'$. Then in \cite{pedicchio}, the commutator $[R,S]$ is defined to be the image under $p''$ of the intersection $\Delta_{R,S}\cap R'$. It is then proved in \cite{pedicchio} that this agrees with the definition of commutators from \cite{fm} at least for Maltsev varieties, so in particular all the properties of the commutators from \cite{fm} hold.

If $[R,S]=\mic_M$ (the smallest congruence diagonal$:M\into M^2$), then the
congruences $R$ and $S$ are said to \emph{centralize} each other. In fact,
$[R,S]$ is the smallest among those congruences $T$ on $M$ for which the
congruences $R/T$ and $S/T$ (on $M/T$) centralize each other. A congruence $R$
on $M$ is called \emph{abelian} if it centralizes itself, i.~e.
$[R,R]=\mic_M$, and \emph{central} if it is centralized by the largest
congruence $\mac_M$ (the identity$:M^2\into M^2$), i.~e. $[\mac_M,R]=\mic_M$.
The \emph{center} $\zeta(M)$ of $M$ is the largest central congruence; it
always exists (in fact more generally for any $R$ and $S$ always exists the
largest congruence $T$ with $[T,R]\le S$). One then defines, generalizing the
usual notions, a central series for a model $M$ to be a chain of congruences
$\mic_M=R_0\le R_1\le...\le R_n=\mac_M$ such that for all $i$ one has
$R_{i+1}/R_i\le\zeta(M/R_i)$ (equivalently, $[\mac_M,R_{i+1}]\le R_i$), i.~e.
$R_{i+1}/R_i$ is a central congruence on $M/R_i$. A model is called
\emph{abelian} if $\zeta(M)=\mac_M$ (equivalently, $[\mac_M,\mac_M]=\mic_M$)
and $n$ \emph{stage nilpotent}, or just $n$-nilpotent, if it has a central
series of length $n$. This happens iff either the \emph{upper central series}
$\mic_M=\zeta^0(M)\le\zeta^1(M)\le\zeta^2(M)\le...$ ends with
$\zeta^n(M)=\mac_M$ or the \emph{lower central series}
$\mac_M=\Gamma^0(M)\ge\Gamma^1(M)\ge\Gamma^2(M)\ge...$ ends with
$\Gamma^n(M)=\mic_M$; here, $\zeta^{n+1}(M)/\zeta^n(M)=\zeta(M/\zeta^n(M))$ 
and $\Gamma^{n+1}(M)=[\mac_M,\Gamma^n(M)]$. Just as in the case of groups,
algebras, etc., the $\Gamma^n$ are functorial (and the $\zeta^n$ are not). A
theory is called $n$-nilpotent (abelian for $n=1$) if all of its models are. 

Let us give another equivalent construction of the commutator in Maltsev
varieties. 

\prp{comm}
{\sl For any congruences $R$, $S$ on an object $A$ in a Maltsev variety,
there is a pushout square
\begin{equation*}
\begin{split}
\xyma{
R\sqcup_AS\ar[r]^{p_{R,S}}
\ar[d]^{q_{R,S}}          &A
                           \ar@{->>}[d]\\
R\sqcap_AS  \ar[r]^{m_{R,S}}  &A/[R,S]\ ,
}
\end{split}\tag{\dag}
\end{equation*}
where $p_{R,S}:R\sqcup_AS\to A$ is induced by the pair
$((x,y)\mapsto x):R\to A$, $((x,y)\mapsto y):S\to A$ and
$q_{R,S}:R\sqcup_AS\to R\sqcap_AS$ is induced by the pair
$((x,y)\mapsto((x,y),(y,y))):R\to R\sqcap_AS$, $((x,y)\mapsto((x,x),(x,y))):S\to R\sqcap_AS$.}
\begin{proof}
Observe that $R\sqcap_AS$ consists of elements of the form $((x,y),(y,z))$
with $(x,y)\in R$ and $(y,z)\in S$. Let $m$ be any Maltsev operation in our
variety, then
\begin{align*}
((x,y),(y,z))&=(m((x,y),(y,y),(y,y)),m((y,y),(y,y),(y,z)))\\
             &=m(((x,y),(y,y)),((y,y),(y,y)),((y,y),(y,z)))\\
             &=m(q_{R,S}i_R(x,y),q_{R,S}i(y,y),q_{R,S}i_S(y,z))\\
             &=q_{R,S}m(i_R(x,y),i(y,y),i_S(y,z)),
\end{align*}
where $i_R$, $i_S$ are the canonical coproduct inclusions and $i$ stands for
any of them.

This shows first of all that $q_{R,S}$ is surjective, so if one forms a
pushout square as above, one gets for the right vertical map the quotient
$A\onto A/T$ for some congruence $T$ on $A$. Moreover it follows that the
induced homomorphism $R\sqcap_AS\to A/T$ in this pushout maps any element
$((x,y),(y,z))=q_{R,S}m(i_R(x,y),i(y,y),i_S(y,z))$ to the $T$-equivalence
class of the element $p_{R,S}m(i_R(x,y),i(y,y),i_S(y,z))$, which equals
$m(p_{R,S}i_R(x,y),p_{R,S}i(y,y),p_{R,S}i_S(y,z))$ $=$ $m(x,y,z)$. 

It follows that one has a commutative diagram
\begin{equation*}
\begin{split}
\xyma{
&R\sqcup_AS\ar[rr]^{p_{R,S}}\ar@{->>}[ldd]
\ar[d]^{q_{R,S}}          &&A\ar[ldd]
                           \ar@{->>}[d]\\
&R\sqcap_AS\ar@{-->}[urr]_m \ar[rr]\ar[ldd]  &&A/T\ar@{=}[dl]\\
R/T\sqcup_{A/T}S/T\ar[rr]^{p_{R/T,S/T}}
\ar@{->>}[d]_{q_{R/T,S/T}}      &&A/T\\
R/T\sqcap_{A/T}S/T\ar@{-->}[urr]_m&&,
}
\end{split}\tag{\ddag}
\end{equation*}
where dashed lines denote maps which are not necessarily homomorphisms.
Indeed, we just showed that the triangle
$$
\xyma{
&&&A\ar@{->>}[d]\\
&R\sqcap_AS\ar@{-->}[urr]_m \ar[rr]&&A/T
}
$$
commutes; whereas the parallelogram
$$
\xyma{
&&&A\ar[ldd]\\
&R\sqcap_AS\ar@{-->}[urr]_m\ar[ldd]\\
&&A/T\\
R/T\sqcap_{A/T}S/T\ar@{-->}[urr]_m&&
}
$$
commutes simply because $A\onto A/T$ is a homomorphism.

Now (\ddag) shows that composing the map $m:R/T\sqcap_{A/T}S/T\to A/T$ with a
surjective homomorphism $R\sqcup_AS\onto R/T\sqcap_{A/T}S/T$ produces a
homomorphism $R\sqcup_AS\to A/T$. It then follows that
$m:R/T\sqcap_{A/T}S/T\to A/T$ is a homomorphism too. Thus by \ref{prp:cent}
$R/T$ and $S/T$ centralize each other, i.~e. $[R,S]\subseteq T$.

Conversely, since $R/[R,S]$ and $S/[R,S]$ centralize each other, there is a
homomorphism $m:(R/[R,S])\sqcap_{A/[R,S]}(S/[R,S])\to A/[R,S]$ as in
\ref{prp:cent}. Composing it with the product of quotient maps
$R\sqcap_AS\to(R/[R,S])\sqcap_{A/[R,S]}(S/[R,S])$ gives a homomorphism
$m_{R,S}$, and to say that it fits in a commutative square as $(\dag)$ above
is precisely the same as to say that $m$ satisfies the Maltsev identities.
This shows that there is a homomorphism $A/T\to A/[R,S]$ under $A$, i.~e. that
$T\subseteq[R,S]$. 
\end{proof}
\rem{comm}
Note that this proposition hints at another possibility of the construction of commutators in varieties more general than Maltsev: given congruences $R$ and $S$ on an algebra $A$, a new congruence $[R,S]$ is uniquely determined by requiring existence of a pushout square
\begin{equation*}
\begin{split}
\xyma{
R\sqcup_AS\ar[r]^{p_{R,S}}
\ar@{->>}[d]^{q_{R,S}} &A
                       \ar@{->>}[d]\\
R\square_AS  \ar[r]    &A/[R,S]\ .
}
\end{split}\tag{\dag}
\end{equation*}
Here we used notation from \ref{prp:comm} above and denoted by $R\square_AS$ the image of $q_{R,S}:R\sqcup_AS\to R\sqcap_AS$.

It is not difficult to determine syntactic content of the above definition of the commutator: the above $[R,S]$ is the smallest congruence with the property that for any two operations $u(x_1,...,x_m,y_1,...,y_n)$, $v(x_1,...,x_p,y_1,...,y_q)$ in our variety and any $a_1,...,a_m$, $b_1,...,b_n$, $c_1,...,c_p$ and $d_1,...,d_q$ in $A$ one has $u(a_1,...,a_m,b_1,...,b_n)\sim_{[R,S]}v(c_1,...,c_p,d_1,...,d_q)$ whenever there exist
$a_i'\sim_Ra_i$, $i=1,...,m$, $b_i'\sim_Sb_i$, $i=1,...,n$, $c_i'\sim_Rc_i$, $i=1,...,p$ and $d_i'\sim_Sd_i$, $i=1,...,q$ satisfying the equalities
$u(a_1',...,a_m',b_1,...,b_n)=v(c_1',...,c_p',d_1,...,d_q)$, $u(a_1,...,a_m,b_1',...,b_n')=v(c_1,...,c_p,d_1',...,d_q')$, and
$u(a_1',...,a_m',b_1',...,b_n')=v(c_1',...,c_p',d_1',...,d_q')$.

This kind of condition has been considered by universal algebraists as a ``better'' behaved than the usual one for non-Maltsev varieties. In fact in \cite{quackenbush} a whole infinite sequence of conditions is considered, each stronger than previous, and ours is the second in the row.

Fortunately we are confined to the realm of Maltsev varieties. In fact all we have to know about commutators is the following fact, which is crucial for what follows: 
\prp{cent}
{\sl
Congruences $R$ and $S$ on a model $M$ of a Maltsev theory centralize each
other if and only if there is a homomorphism $m$ from the submodel
$$
R\sqcap_MS\cong\setof{(x,y,z)\in M^3}{(x,y)\in R,(y,z)\in S}
$$
of $M^3$ to $M$ which satisfies $xSm(x,y,z)Rz$ for any $xRySz$ and is Maltsev,
i.~e. $m(x,y,y)=x$ for any $xRy$, and $m(y,y,z)=z$ for any $ySz$. Then
restriction of any Maltsev operation $p$ of the theory to that submodel
coincides with this homomorphism, is associative and commutative. More
precisely, $p(u,v,p(x,y,z))=p(p(u,v,x),y,z)$ holds if $uRv$, $ySz$, and either
$vSx$ or $xRy$, while $p(x,y,z)=p(z,y,x)$ holds for any} $xRySz$. 
\begin{proof}
The first statement can be found in \cite{pedicchio}, see Lemma 2.11 there; we
only prove the second (the proof is essentially the same as in
\cite{johnstone}). 

Take any Maltsev operation $p$. Then homomorphicity of $m$ means
$$
m(p(x_1,x_2,x_3),p(y_1,y_2,y_3),p(z_1,z_2,z_3))
=p(m(x_1,y_1,z_1),m(x_2,y_2,z_2),m(x_3,y_3,z_3)),
$$
for any $x_iRy_iSz_i$ ($i=1,2,3$). Taking here $y_1=z_1=x_2=y_2=z_2=x_3=y_3$
then gives $m(x_1,y_2,z_3)=p(x_1,y_2,z_3)$ for any $x_1Ry_2Sz_3$. Now for any
$uRv$, $xRySz$ one has 
\begin{multline*}
p(u,v,p(x,y,z))=p(m(u,v,v),m(v,v,v),m(x,y,z))\\
=m(p(u,v,x),p(v,v,y),p(v,v,z))=m(p(u,v,x),y,z)=p(p(u,v,x),y,z),
\end{multline*}
and similarly for $uRvSx$, $ySz$. Whereas taking any $xRySz$,
\begin{multline*}
p(z,y,x)=p(m(y,y,z),m(y,y,y),m(x,y,y))\\
=m(p(y,y,x),p(y,y,y),p(z,y,y))=m(x,y,z)=p(x,y,z).
\end{multline*}
\end{proof}

\crl{torcen}
{\sl
A congruence $R$ on a model $M$ of a Maltsev theory $\t$ is abelian
if and only if the morphism $M\to M/R$ is an abelian extension in the sense of
\ref{dfn:abcen}. Furthermore $R$ is central if and only if $M\to M/R$ is (the
trivial case $M=\varnothing$ excluded).}
\begin{proof}
The first statement is immediate by \ref{prp:maltor} and \ref{prp:cent},
and the second by \ref{prp:malcen}.
\end{proof}

In fact, this is a special case of a statement which deals with arbitrary
$R$, $S$ centralizing each other. This requires generalizing herd structures
to s.~c. \emph{herdoids}; see \cite{pedicchio}.

\crl{cenab}
{\sl
A congruence $R$ on a non-empty model $M$ of a Maltsev theory $\t$ is central
if and only if there is an internal abelian group $A$ in $\t\m$, an action of
$A$ on $M$, and an isomorphism $\ph:R\to A\x M$ fitting in the commutative
triangle}
$$
\xyma{
R\ \ar@{>->}[dr]\ar[dd]_\ph^{\cong}\\
&M\x M\\
A\x M\ar[ur]_{\textrm{\ \ \ \ \ (projection, action)}}.
}
$$
\begin{proof}
This follows easily from \ref{crl:torcen}, in view of \ref{prp:malcen}.
\end{proof}

Since linear, respectively untwisted, extensions of a theory $\t$ by
cartesian natural systems (resp. bifunctors) are, by \ref{prp:lin=tor},
none other than abelian, resp. central, objects of $\Th/\t$, one might
expect that they are related to abelian, resp. central, extensions in
$\t\m$. For Maltsev theories, a link between these notions is provided by

\thm{th=mod}
{\sl
For a morphism $P:\t'\to\t$ of Maltsev theories, the following conditions are
equivalent: 
\begin{itemize}
\item[i)]
$P$ is a linear (respectively, untwisted) extension;
\item[ii)] for all $\t'$-models $M$ the homomorphisms $\eta_M:M\to M_P$ are
abelian (respectively, central) extensions in $\t'\m$ and, moreover, the following condition is satisfied: 
\begin{itemize}
\item[$\diamondsuit$]
for any Maltsev operation $p$ in $\t'$, any $u,v:\t'(n)\to\t'(k)$ in $\t'\m$ with
$P(u)=P(v)$, and any $x,y\in\t'(n)$ with $\eta_{\t'(n)}(x)=\eta_{\t'(n)}(y)$, one has
$$
p(u(x),v(x),v(y))=u(y).
$$
\end{itemize}
\end{itemize}
}

\begin{proof}

\ 

i) $\then$ ii):

By \ref{prp:lin=tor}, the morphism $P$ is a linear extension if and only if it
is full and there is a functor $m:\t'\x_\t\t'\x_\t\t'\to\t'$ over $\t'$ which
is a commutative associative Maltsev operation in $\Th/\t$. Identifying $\t'$
with the opposite of the category of finitely generated free models via the Yoneda
embedding, the action of $m$ on $\hom(X^n,X)$ may be viewed as a commutative
associative Maltsev operations
$m_n:\t'(n)\x_{\t(n)}\t'(n)\x_{\t(n)}\t'(n)\to\t'(n)$ over $\t(n)$. Now
functoriality of $m$ means under the above identification that for any
$u,v,w\in\hom_{\t\m}(\t'(n),\t'(k))\approx\t'(k)^n$ with $P(u)=P(v)=P(w)$, and
any $x,y,z\in\hom_{\t\m}(\t'(i),\t'(n))\approx\t'(n)^i$ with $P(x)=P(y)=P(z)$
one has 
\begin{equation}
m_k^i(ux,vy,wz)=m_k^n(u,v,w)m_n^i(x,y,z).\tag{$\circ$}
\end{equation}
This then shows that each $m_n$ is a homomorphism: taking $i=1$ and $x=y=z$ in
$(\circ)$ gives 
$$
m_k(ux,vx,wx)=m_k^n(u,v,w)m_n(x,x,x)=m_k^n(u,v,w)x
$$
which, in terms of $\t'$ again, means
$$
m_k(x(u_1,...,u_k),x(v_1,...,v_k),x(w_1,...,w_k))
=x(m_k(u_1,v_1,w_1),...,m_k(u_n,v_n,w_n)).
$$
So $\eta_{\t'(n)}$ is a linear extension, and $m_n$ coincides with the
restriction of any Maltsev operation on $\t'(n)$. If moreover $P$ is
untwisted, then, as in \ref{prp:malcen}, this $m_n$ is defined on
$\t'(n)\x_{\t(n)}\t'(n)\x\t'(n)$, and $\eta_{\t'(n)}$ is central. 

Now taking $i=1$, $v=w$ and $x=y$ in ($\circ$) gives
$$
m_k(u(x),v(x),v(z))=m_k^n(u,v,v)m_n(x,x,z)=u(z),
$$
which, since the $m$'s coincide with the restrictions of Maltsev operations,
gives $\diamondsuit$. 

Now for a general free model $\t'(S)$, the homomorphism $\eta_{\t'(S)}$ is a
colimit of a filtered diagram of those of the form $\eta_{\t'(n)}$, just as in
the proof of \ref{prp:subvar}. Since filtered colimits commute with finite
limits and are created by the forgetful functors, it follows that the
collection of the $m_n$ on the $\t'(n)$ over $\t(n)$ give rise to one $m_S$ on
$\t'(S)$ over $\t(S)$, so $\eta_{\t'(S)}$ is abelian, resp. central, whenever
all the $\eta_{\t'(n)}$ are. 
 
Finally, consider any $\t'$-model $M$. Let us choose a surjective homomorphism
$q:\t'(S)\onto M$, so that $M=\t'(S)/R_M$ for some congruence $R_M$ on
$\t'(S)$. By \ref{prp:subvar}, both $\eta_{\t'(S)}$ and $\eta_M$ are
surjective, so also $\t(S)=\t'(S)/R_\t$, $M_P=\t'(S)/R$, for certain
congruences $R_\t$, $R$, with $R_M\subseteq R$. In fact, the pushout condition
from \ref{prp:subvar} shows that $R=R_M\lor R_\t$ in the lattice of
congruences on $\t'(S)$. And since $\t'$ is Maltsev, in fact $R=R_M\o R_\t$.
Thus for any $u_1,u_2,\dots\in\t'(S)$, one has $u_1Ru_2R\dots$ iff
$u_1R_Mv_1R_\t u_2R_Mv_2R_\t\dots$, for some $v_1,v_2,\dots$. We then conclude
that for any $x_1,x_2,\dots\in M$, one has $\eta_M(x_1)=\eta_M(x_2)=\dots$ iff
there are $u_i\in\t'(S)$ with $x_i=q(u_i)$ $(i=1,2,\dots)$ and 
$\eta_{\t'(S)}(u_1)=\eta_{\t'(S)}(u_2)=\dots$. 

We then define the Maltsev operation $m_M:M\x_{M_P}M\x_{M_P}M\to M$
over $M_P$ (respectively $M\x_{M_P}M\x M\to M$ for untwisted $P$) by
$m_M(x_1,x_2,x_3)=qm_S(u_1,u_2,u_3)$, for some
$(u_1,u_2,u_3)$ in $\t'(S)\x_{\t(S)}\t'(S)\x_{\t(S)}\t'(S)$
(respectively, in $\t'(S)\x_{\t(S)}\t'(S)\x\t'(S)$) with $q(u_i)=x_i$ -- which
exist by the preceding argument. This is legitimate since for any other choice
$v_i$ one would have $v_iR_Mu_i$, hence $m_S(v_1,v_2,v_3)R_Mm_S(u_1,u_2,u_3)$.
This since $R_M$ is a submodel of $\t'(S)^2$, while $m_S$, by \ref{prp:cent},
coincides with the restriction of any Maltsev operation. Homomorphicity,
Maltsev identities, associativity and commutativity of $m_M$ now follow from
those of $m_S$. 

\

ii) $\then$ i):

By \ref{prp:maltor}, \ref{prp:lin=tor}, and \ref{prp:simtor}, to prove that
$P$ is a linear (resp., untwisted) extension, it suffices to construct a
family of commutative associative Maltsev operations, denoted
$(u,v,w)\mapsto u-v+w$, from 
$$
\hom_{\t'}(X^n,X^k)\x_{\hom_\t(X^n,X^k)}\hom_{\t'}(X^n,X^k)\x_{\hom_\t(X^n,X^k)}\hom_{\t'}(X^n,X^k)
$$
--- respectively, from
$$
\hom_{\t'}(X^n,X^k)\x_{\hom_\t(X^n,X^k)}\hom_{\t'}(X^n,X^k)\x\hom_{\t'}(X^n,X^k)
$$
--- to $\hom_{\t'}(X^n,X^k)$, which define a functor $P\x P\x P\to P$ over
$\t$ in $\Th/\t$, (resp., and also satisfy the conditions as in
\ref{prp:simtor}). For that, choose some Maltsev operation
$p\in\hom_{\t'}(X^3,X)$, and put
$$
\bar u-\bar v+\bar w=(p(u_1,v_1,w_1),...,p(u_k,v_k,w_k)),
$$
for any $\bar u=(u_1,...,u_k),\bar v=(v_1,...,v_k),\bar w=(w_1,...,w_k)\in\hom_{\t'}(X^n,X^k)$.
Functoriality then amounts to
$$
p(u\bar u,v\bar v,w\bar w)=p(up(\bar u,\bar v,\bar w),vp(\bar u,\bar v,\bar w),wp(\bar u,\bar v,\bar w)),
$$
for $u,v,w\in\hom_{\t'}(X^k,X)$ and $\bar u,\bar v,\bar w$ as above,
whenever $P(\bar u)=P(\bar v)=P(\bar w)$ and $P(u)=P(v)=P(w)$.
Whereas the conditions from \ref{prp:simtor} become
$$
u(p(u_1,v_1,w_1),...,p(u_k,v_k,w_k))=p(u\bar u,u\bar v,u\bar w)\textrm{ and }
p(u,v,w)\bar u=p(u\bar u,v\bar u,w\bar u),
$$
if $P(\bar u)=P(\bar v)$ and $P(u)=P(v)$.

Let us use the Yoneda embedding to identify $\hom_{\t'}(X^i,X^j)$ with
$\t'(i)^j$. Then, what we have to prove is this: for any $u,v,w\in\t'(k)$ and
any homomorphisms $\bar u,\bar v,\bar w:\t'(k)\to\t'(n)$, the equality 
$$
p(\bar u(u),\bar v(v),\bar w(w))=p(\bar u(p(u,v,w)),\bar v(p(u,v,w)),\bar w(p(u,v,w)))
$$
holds whenever $P(\bar u)=P(\bar v)=P(\bar w)$ and
$\eta_{\t'(n)}(u)=\eta_{\t'(n)}(v)=\eta_{\t'(n)}(w)$ (resp., also
\begin{gather*}
p(\bar u,\bar v,\bar w)(u)=p(\bar u(u),\bar v(u),\bar w(u)),\\
\bar u(p(u,v,w))=p(\bar u(u),\bar u(v),\bar u(w)),
\end{gather*}
when $P(\bar u)=P(\bar v)$ and $\eta_{\t'(n)}(u)=\eta_{\t'(n)}(v)$).

The functoriality condition, using that $\bar u$, $\bar v$, $\bar w$ are
homomorphisms and hence commute with the operation $p$, is equivalent to 
$$
p(\bar u(u),\bar v(v),\bar w(w))=
p(p(\bar u(u),\bar u(v),\bar u(w)),p(\bar v(u),\bar v(v),\bar v(w)),p(\bar w(u),\bar w(v),\bar w(w))).
$$
Now recall that for the given elements $p$ is a commutative associative
Maltsev homomorphism, so, switching to additive notation and using
\ref{lem:asmal}, the equality in question becomes 
$$
\bar u(u)-\bar v(v)+\bar w(w)=
\bar u(u)-\bar u(v)+\bar u(w)-\bar v(w)+\bar v(v)-\bar v(u)+\bar w(u)-\bar w(v)+\bar w(w).
$$
Let us replace $\bar u(u)-\bar v(v)+\bar w(w)$ by
$\bar u(u)-\bar v(v)+\bar v(v)-\bar v(v)+\bar w(w)$ and then substitute, using
$\diamondsuit$, 
$$
-\bar v(v)=-(\bar v(w)-\bar u(w)+\bar u(v))
$$
in the first occurrence and
$$
-\bar v(v)=-(\bar w(v)-\bar w(u)+\bar v(u))
$$
in the second. We then obtain
$$\setlength\arraycolsep{0pt}
\begin{array}{rrrrrrrr}
 \bar u(u)&-(\bar v(w)-\bar u(w)& +\bar u(v)&)+\bar v(v)&-(\bar w(v)& -\bar w(u)+\bar v(u)&)+\bar w(w)&\\
=\bar u(u)& -\bar u(v)+\bar u(w)& -\bar v(w)& +\bar v(v)& -\bar v(u)& +\bar w(u)-\bar w(v)& +\bar w(w)&\\
=\bar u(u)& -\bar u(v)+\bar u(w)&-(\bar v(u)& -\bar v(v)& +\bar v(w)&)+\bar w(u)-\bar w(v)& +\bar w(w)&,
\end{array}
$$
as required.

As for the last two equalities, the first is trivial, and the second follows
since $\bar u$ is a homomorphism. 
\end{proof}

Outside the realm of Maltsev theories there exist linear extensions
$P:\t'\to\t$ such that not all the $\eta$'s are abelian. 

\subsection*{Example}
Let $M$ be the multiplicative monoid $\set{1,0}$, i.~e. $1\cdot1=1$,
$1\cdot0=0\cdot1=0\cdot0=0$. Consider the natural system $D$ on it, in the
sense explained in \ref{exm}.\ref{mon}, given by $D_1={0}$,
$D_0=\z/2\z\os\z/2\z$, with $0(x,y)=(y,y)$ and $(x,y)0=(0,0)$ for $(x,y)\in
D_0$. Then the trivial linear extension of $M$ by $D$ is the monoid 
$M\ax D=\set{(1,0)}\cup\set{(0,(0,0)),(0,(0,1)),(0,(1,0)),(0,(1,1))}$, with
$(1,0)\cdot(0,(x,y))=(0,(x,y))\cdot(1,0)=(0,(x,y))$ and
$(0,(x,y))\cdot(0,(x',y'))=(0,(x,y)0+0(x',y'))=(0,(y',y'))$. 
For brevity, let us redenote this as $M\ax D=\set{1,00,10,01,11}$, so that $1$
is the unit and 
\begin{align*}
&00\cdot00=10\cdot00=01\cdot00=11\cdot00=00\cdot10=10\cdot10=01\cdot10=11\cdot10=00,\\
&00\cdot01=10\cdot01=01\cdot01=11\cdot01=00\cdot11=10\cdot11=01\cdot11=11\cdot11=11.
\end{align*}

As in \ref{exm}.\ref{mon}, $D$ extends uniquely to a natural system on $\t_M$
in such a way that the morphism $P:\t_{M\ax D}\to\t_M$, induced by the
projection $M\ax D\to M$, is the trivial linear extension by that system. Now
let $S=\set{1,*0,01,11}=M\ax D/(00\sim10)$ be the $M\ax D$-set obtained by
identifying $00$ and $10$ in $M\ax D$, acting on itself from the left via
$\cdot$; then clearly $S_P=(M\ax D)_P=M$, 
with $\eta_S(1)=1$ and $\eta_S(*0)=\eta_S(01)=\eta_S(11)=0$. Suppose $\eta_S$
is abelian. Then there exists a Maltsev operation $m:S\x_MS\x_MS\to S$ over
$M$. Since it must be a morphism of $M\ax D$-sets, one must have in particular
$10\cdot m(*0,01,11)=m(10\cdot*0,10\cdot01,10\cdot11)=m(*0,11,11)=*0$. This is
only possible if $m(*0,01,11)=*0$. Then $m$ cannot be associative, since this
would imply $11=m(11,*0,*0)=m(11,m(*0,01,11),*0)=m(11,11,m(01,*0,*0))=01$. 

\ 

On the other hand, if all the $\eta$'s are abelian, the morphism of
(non-Maltsev) theories can still fail to be a linear extension -- consider the
morphism $\s.\to\1$ from the theory of pointed sets to the terminal theory.
Clearly for any non-empty set $S$ the map $S\to1$ is abelian, as $S$ can be
equipped with an abelian group structure. But $\s.$ itself clearly cannot have
any functorial Maltsev operation $m:\s.\x_\1\s.\x_\1\s.\to\s.$ since
$m(f,g,h)m(f',g',h')=m(ff',gg',hh')$ implies that
$ee'=m(e,1,1)m(1,1,e')=m(e,1,e')=m(1,1,e')m(e,1,1)=e'e$ for any endomorphisms
of any object, while $\s.$ has noncommutative endomorphism monoids. 

\crl{clem}
{\sl
Let $P:\t'\to\t$ be a morphism of Maltsev theories such that for each $M$ in
$\t'\m$ the morphism $\eta_M:M\to M_P$ is a central extension and moreover the
subvariety $U_P(\t\m)$ of $\t'\m$ contains all abelian $\t'$-models. Then $P$
is an untwisted linear extension.}
\begin{proof}
In view of the previous theorem we just have to check the condition
$\diamondsuit$. It is true in more generality, in fact: let $u,v:M'\to M$ be
any homomorphisms in $\t'\m$ with $P(u)=P(v)$, and let $x,y\in M'$ be any
elements with $\eta_{M'}(x)=\eta_{M'}(y)$. We thus have $(u,v):M'\to
R_M=M\x_{M_P}M$. Then by \ref{crl:cenab} we know that there is an action
$A(M)\x M\to M$ of an internal abelian group $A(M)$ on $M$ and an isomorphism
$R_M\cong A(M)\x M$. Thus $(u,v)$ give rise to a homomorphism
$u-v:M'\to A(M)$. Now obviously any internal abelian group in $\t'\m$ is an
abelian model of $\t'$, so by hypothesis $A(M)$ belongs to the image of $U_P$;
hence $u-v$ factors through $\eta_{M'}$ and we obtain $(u-v)(x)=(u-v)(y)$.
Using now the action of $A(M)$, \ref{crl:cenab} gives
$u(x)-v(x)+z=u(y)-v(y)+z$ for any $z$ in $M$, in particular taking $z=v(y)$
gives $\diamondsuit$. 
\end{proof}

\subsection*{Remark}
It is easy to identify the corresponding bifunctor on $\t$: one can show that
it is given by 
$$
D(X^n,X^k)=\hom_{\t\m}(\t(k),A(\t'(n)))\cong A(\t'(n))^k.
$$

\crl{nil}
{\sl
A theory $\t$ is $n$-nilpotent Maltsev if and only if there is a tower of
untwisted linear extensions of theories
$$
\t=\t_n\to\t_{n-1}\to...\to\t_1
$$
where $\t_1$ is an abelian Maltsev theory.}
\begin{proof}
Fix a theory $\t$ and some $n$, and let $\t_{(n)}$ be the theory with
$$
\hom_{\t_{(n)}}(X^k,X^l)=\hom_{\t\m}(\t(l)/\Gamma^n(\t(l)),\t(k)/\Gamma^n(\t(k))).
$$
It is easy to see that there is a bijection
$$
\hom_{\t\m}(M,N/\Gamma^n(N))\approx\hom_{\t\m}(M/\Gamma^n(M),N/\Gamma^n(N))
$$
for any models $M$, $N$, hence it follows that $\t_{(n)}$ is an $n$ stage
nilpotent theory and that sending a model $M$ to $M/\Gamma^n(M)$ restricts to
a morphism of theories $P:\t\to\t_{(n)}$. Moreover all models in the image of
$U_P$ are $n$-nilpotent, and for any $M$ the morphism $\eta_M:M\to M_P$ is the
quotient $q_n(M):M\to M/\Gamma^n(M)$.

Now suppose $\t$ is itself $(n+1)$-nilpotent and Maltsev; then, the above
construction gives a morphism of theories $\t\to\t_n$ whose structure maps
$\hom_\t(X^k,X^l)$ $\to$ $\hom_{\t_n}(X^k,X^l)$ coincide with
$q_n(\t(k))^l:\t(k)^l\onto(\t(k)/\Gamma^n(\t(k)))^l$, with $q_n$ as above.
The fact that $\t(k)$ is $(n+1)$-nilpotent means according to \ref{crl:torcen}
that $q_n$ is a central extension in the sense of \ref{dfn:abcen}. Thus,
\ref{crl:clem} shows that $\t\to\t_n$ is an untwisted linear extension. By
induction, one thus gets the ``only if'' part. 

For the ``if'' part, suppose we are given an untwisted linear extension
$P:\t\to\t_n$ of theories, and $\t_n$ is $n$-nilpotent and Maltsev. Then by
\ref{prp:mal} $\t$ is also Maltsev. Furthermore by \ref{thm:th=mod}, for any
$\t$-model $M$ the homomorphism $\eta_M:M\to M_P$ is a central extension,
hence a quotient by a central congruence, by \ref{crl:torcen}. Since $P(M)$ is
$n$-nilpotent, it follows that $M$ is $(n+1)$-nilpotent.
\end{proof}

\subsection*{Remark} Let us also note in this connection that any nilpotent
(in fact, also \emph{solvable}) variety with modular congruence lattices is
Maltsev -- see 10.1 in \cite{gumm}. 

\section{Abelian theories}

We finish with some information on the structure of abelian Maltsev theories
and linear extensions between them. We will first deal with abelian theories
without constants, i.~e. nullary operations. Abelian theories with constants
are much simpler and will be treated in the end.

It follows from \ref{prp:cent} that a model $A$ of a Maltsev theory is abelian
iff there is a homomorphism $m:A^3\to A$ satisfying the Maltsev identities,
and then any Maltsev operation on $A$ coincides with $m$, is commutative and
associative. It will be convenient to fix such an operation throughout and
denote it by $m(x,y,z)=x+_yz$. Thus all models of an abelian theory are
\emph{abelian herds}, and all of their operations are homomorphisms of abelian
herds. To describe such theories, we need the following notions. 

\dfn{cross}
A \emph{left linear form} consists of an associative ring $R$ with unit,
a left $R$-module $M$, and a homomorphism $\d:M\to R$ of left $R$-modules.

In fact usually we will omit the word ``left'', as it is customary with
modules. 

\dfn{aff}
An \emph{affinity} over a linear form $\d:M\to R$ is an abelian herd
$A$ together with maps $R\x A\x A\to A$ and $M\x A\to A$, denoted,
respectively, $(r,a,b)\mapsto r_ab$ and $(x,a)\mapsto\ph_a(x)$, such that the
following identities hold:
\begin{itemize}
\item
For each $a\in A$, the operations $(-)+_a(-)$ and $(-)_a(-)$ turn $A$ into a
left $R$-module (with zero $a$), and $\ph_a$ into a module homomorphism. In
other words, for any $a,b,c,d\in A$, $r,s\in R$, $x,y\in M$ one has
\begin{align*}
&b+_a(c+_ad)=(b+_ac)+_ad,\\
&a+_ab=b,\\
&b+_ac=c+_ab,\\
&b-_ab=a,\\
&r_a(b+_ac)=r_ab+_ar_ac,\\
&(r+s)_ab=r_ab+_as_ab,\\
&1_ab=b,\\
&r_a(s_ab)=(rs)_ab,\\
&\ph_a(x+y)=\ph_a(x)+_a\ph_a(y),\\
&\ph_a(rx)=r_a\ph_a(x),
\end{align*}
where we have denoted $b-_ac=b+_a((-1)_ac)$.
\item (``Coordinate change''.) These structures are related by the identities
\begin{align*}
&b+_{a'}c=((b-_aa')+_a(c-_aa'))+_aa',\\
&r_{a'}b=r_a(b-_aa')+_aa',\\
&\ph_{a'}(x)=\ph_a(x)+_a(1-\d x)_aa'. \tag{*}
\end{align*}
\end{itemize}
A homomorphism between affinities $A$, $A'$ is a map $f:A\to A'$ preserving
all this, i.~e. satisfying
\begin{align*}
&f(a+_bc)=f(a)+_{f(b)}f(c),\\
&f(r_ab)=r_{f(a)}f(b),\\
&f(\ph_a(x))=\ph_{f(a)}(x).
\end{align*}

Obviously the category $\d\aff$ of affinities over a linear form $\d$
is the category of models of a suitable abelian theory $\t_\d$. Here is an
explicit description of this theory.

\prp{crab}
{\sl
The theory $\t_\d$ of affinities over $\d:M\to R$ can be described as follows:
$$
\hom_{\t_\d}(X^n,X)=
\begin{cases}
\varnothing,&n=0;\\
M\x R^{n-1},&n>0.
\end{cases}
$$
The projections $a_0,a_1,a_2,...:X^n\to X$ are given, respectively, by the
elements $\brk{0,0,0,...}$, $\brk{0,1,0,...}$, $\brk{0,0,1,...}$, ...; and,
composition is given by
$$
\brk{x,r_1,r_2,...}(\brk{x_0,s_0,t_0,...},\brk{x_1,s_1,t_1,...},\brk{x_2,s_2,t_2,...},...)=
\brk{x',s',t',...},
$$
where
\begin{align*}
&x'=x+(1-\d(x))x_0+r_1(x_1-x_0)+r_2(x_2-x_0)+...,\\
&s'=(1-\d(x))s_0+r_1(s_1-s_0)+r_2(s_2-s_0)+...,\\
&t'=(1-\d(x))t_0+r_1(t_1-t_0)+r_2(t_2-t_0)+...,\\
...
\end{align*}}
\begin{proof}
Take as basic operations the ternary $(-)+_{(-)}(-)$, the family of binaries
$r_{(-)}(-)$ indexed by $r\in R$, and unaries $\ph_{(-)}(x)$ indexed by $x\in
M$. Using the above identities, one can write any composite of these
operations in the form
$$
\brk{x,r,s,...}(a,b,c,...)
=\ph_a(x)+_ar_ab+_as_ac+_a\cdots
$$
in a unique way. The rest is straightforward verification.
\end{proof}

Define now a morphism of left linear forms from $\d:M\to R$ to
$\d':M'\to R'$ to be an equivariant homomorphism, i.~e. a pair
$(f:R\to R',g:M\to M')$ of additive maps such that the obvious square
commutes, that $f$ is a unital ring homomorphism, and that $g(rx)=f(r)g(x)$
holds for any $r\in R$, $x\in M$. This clearly defines the category $\Lf$ of
left linear forms. We then have 

\thm{crab}
{\sl
The category of abelian Maltsev theories without nullary operations is
equivalent to the category $\Lf$ of left linear forms.}
\begin{proof}
Define the functor $\t_{(-)}:\Lf\to\Th$ by sending an object $\d$ of $\Lf$
to the corresponding theory $\t_\d$ described above in \ref{prp:crab}. It is
clear from that description that any morphism in $\Lf$ determines a
morphism of the corresponding theories.

Conversely, given an abelian Maltsev theory $\t$, define the left linear form
$\d_\t:M_\t\to R_\t$ as follows: let $M_\t$ be the set of all unary operations
of $\t$, with the abelian group structure given by $x+y=m(x,\id,y)$, where $m$
is the Maltsev operation of $\t$. Let $R_\t$ be the set of \emph{convex}
binary operations of $\t$, i.~e. those binary operations $r$ satisfying the
identity $r(a,a)=a$. Define the ring structure on it by taking zero $0$ to be
the first projection, unit $1$ the second projection, addition to be
$r+s=m(r,0,s)$, additive inverse $-r=m(0,r,0)$, and multiplication
$(rs)(a,b)=r(a,s(a,b))$. Let $R_\t$ act on $M_\t$ via $(rx)(a)=r(a,x(a))$, and
let the crossing $M_\t\to R_\t$ be $(\d_\t x)(a,b)=m(x(a),x(b),b)$. It is then
straightforward to check that this defines a left linear form, that any
morphism of theories gives rise to a morphism in $\Lf$ in a functorial way,
and that if one starts from a theory of the form $\t_\d$, then one recovers
the original $\d$ back. Finally for the second way round, observe that for any
operation $u:X^n\to X$ in an abelian theory, with $n>0$, one has 
\begin{align*}
u(a,b,c,...)&=u(a+_aa+_aa+_a\cdots,a+_ab+_aa+_a\cdots,a+_aa+_ac+_a\cdots,...)\\
            &=u(a,a,a,...)+_{u(a,a,a,...)}u(a,b,a,...)+_{u(a,a,a,...)}u(a,a,c,...)+_{u(a,a,a,...)}\cdots\\
            &=u(a,a,a,...)+_a(a+_{u(a,a,a,...)}u(a,b,a,...))+_a(a+_{u(a,a,a,...)}u(a,a,c,...))+_a\cdots\\
            &=x(a)+_ar(a,b)+_as(a,c)+_a\cdots,
\end{align*}
with $x$ in $M_\t$ and $r,s,...$ in $R_\t$. This implies easily that including
$M_\t$ and $R_\t$ in $\t$ extends to an isomorphism of theories from
$\t_{\d_\t}$ to $\t$.
\end{proof}

\subsection*{Remark}
Construction of the ring $R_\t$ from an abelian Maltsev theory $\t$ is
obviously well known to universal algebraists, in a slightly different context
-- see e.~g. \cite{csakany}. It is in fact closely related to the classical
coordinatization construction for geometries. The reader might consult, e.~g.
\cite{gumm} or \cite{fm} for that.

\

Using our description, we can now find out what kind of linear extensions
exist between abelian theories. Indeed, since the Maltsev operation in abelian
theories is unique, they are clearly closed under arbitrary finite limits,
hence \ref{thm:crab} together with \ref{prp:lin=tor} implies that abelian
linear extensions of an abelian theory $\t$ can be identified with torsors
under internal abelian groups in $\Lf/\d_\t$. Thus we just have to describe
torsors under a linear form $\d:M\to R$. Consider one such, given by
\begin{equation}
\tag{E}
\begin{split}
\xyma{
K\ar@{>->}[r]^j\ar[d]^\gd &N\ar@{->>}[r]^q\ar[d]^{\d'} &M\ar[d]^\d\\
B\ar@{>->}[r]^i           &S\ar@{->>}[r]^p             &R.
}
\end{split}
\end{equation}
Now by \ref{prp:maltor} we know that this torsor is equipped with a herd
structure in $\Lf/\d$. Thus we have a Maltsev homomorphism
$m:\d'\x_\d\d'\x_\d\d'\to\d'$ over $\d$; then, by an argument just as in
\ref{prp:cent}, both in $N$ and $S$ one has 
$$
m(x,y,z)=m(x-y+y,y-y+y,y-y+z)=m(x,y,y)-m(y,y,y)+m(y,y,z)=x-y+z.
$$
Thus it follows that $(p,q)$ above is a torsor iff this map is a homomorphism.
One sees easily that this happens iff $B^2=BK=0$. In such case, $B$ becomes
naturally an $R$-$R$-bimodule, $K$ a left $R$-module, and restriction $\gd$ of
$\d'$ to it -- a module homomorphism, via $rb=sb$, $br=bs$, $rk=sk$, for any 
$b\in B$, $k\in K$, $r\in R$ and $s\in S$ with $p(s)=r$. Moreover there is an
$R$-module homomorphism $B\ox_RM\to K$, denoted $(b,m)\mapsto b\cdot m$, given
by $b\cdot m=bn$ for any $n\in N$ with $q(n)=m$. It clearly satisfies
$\gd(b\cdot m)=b\d(m)$. On the whole, one gets a structure which can be
described by 

\dfn{bim}
For a left linear form $\d:M\to R$, a $\d$-\emph{bimodule} consists of an
$R$-$R$-bimodule $B$, a left $R$-module $K$, and $R$-linear maps $\gd:K\to B$
and $\cdot:B\ox_RM\to K$ satisfying $\gd(b\cdot m)=b\d m$ for any $b\in B$,
$m\in M$. It will be denoted $\gd^\cdot=(B\ox_RM\to K\to B)$. 

\ 

Examples of such $\d$-bimodules include $R\ox_RM\cong M\to R$, i.~e. $\d$
itself, $(B\ox_RM\to B\ox_RR\cong B=B)$, for any $R$-$R$-bimodule $B$, which
we denote $\co(B)$, and $0\to K\to0$ for any left $R$-module $K$, which we
denote $K[1]$. 

It is easy to show that also conversely, internal groups in $\Lf/\d$ are
determined by $\d$-bimodules $\gd^\cdot$ as above, and that as soon as the map
$(x,y,z)\mapsto x-y+z$ is a homomorphism from $\d'\x_\d\d'\x_d\d'$ to $\d'$
over $\d$, then there is a torsor structure on $\d'$ under the corresponding
group. We summarize this as follows: 

\prp{crext}
{\sl
For a left linear form $\d:M\to R$, internal groups in $\Lf/\d$ are in
one-to-one correspondence with $\d$-bimodules. Moreover the underlying linear
form of the corresponding group is $\gd\os\d:K\os M\to B\os R$, with
multiplicative structure $(b,r)(b',r')=(br'+rb',rr')$,
$(b,r)(k,m)=(b\cdot m+rk,rm)$. Torsors under this group are in one-to-one
correspondence with diagrams such as $(\mathrm E)$ above, where
$B\into S\onto R$ is a singular extension, i.~e. the ideal $i(B)$ has zero
multiplication in $S$ and the induced $R$-$R$-bimodule structure coincides
with the original one, and 
moreover $i(B)j(K)=0$, the induced $R$-module structure on $K$ is the original
one, i.~e. $j(p(s)k)=sj(k)$, and finally, the induced action $B\ox_RM\to K$
coincides with the original one, i.~e. $j(b\cdot q(n))=i(b)n$.
}\qed

\ 

Translating now all of the above from $\Lf$ to abelian theories, we conclude

\prp{abext}
{\sl
For a left linear form $\d:M\to R$, each internal group $\A=(\gd\os\d\to\d)$
in $\Lf/\d$ corresponding to the $\d$-bimodule $\gd^\cdot=(B\ox_RM\to K\to B)$
as above, gives rise to a natural system $D^\A$ on the corresponding abelian
theory $\t_\d$. Explicitly, one has 
$$
D^\A_{\brk{x,r_1,...,r_{n-1}}}=K\os B^{n-1},
$$
with actions given by restricting those in \ref{prp:crab} for $\t_{\gd\os\d}$
to $K\os B^{n-1}\subseteq(K\os M)\x(B\os R)^{n-1}$.
 
A natural system on $\t_\d$ is of this form iff all linear extensions by it
are again abelian theories. 
} \qed

\ 

In view of this, we will in what follows identify internal groups $\A$ in
$\Lf/\d$ with $\d$-bimodules and with the corresponding natural systems $D^\A$
on $\t_\d$. In particular, equivalence classes of extensions of $\t_\d$ by
$D^\A$ form, by \cite{jp}, an abelian group isomorphic to $H^2(\t_\d;D^\A)$,
which we can as well denote $H^2(M\to R;\ B\ox_RM\to K\to B)$, or just by
$H^2(\d;\gd^\cdot)$.

Now from \cite{jp} we know that any short exact sequence
$\gd'\into\gd^\cdot\onto\gd''$ induces the exact sequence 
$$
0\to H^0(\d;\gd')\to...\to H^1(\d;\gd'')\to H^2(\d;\gd')\to H^2(\d;\gd^\cdot)\to H^2(\d;\gd'')\to...
$$
which one can use to reduce investigation of cohomologies, in particular
linear extensions by a $\d$-bimodule, to those by more ``elementary'' ones. In
particular, observing the diagrams 
$$
\xyma{
\ker(\gd)\ar@{>->}[r]\ar[d] &K\ar@{->>}[r]\ar[d]^\gd &\im(\gd)\ar@{>->}[d]^\iota\\
0\ar[r]                     &B\ar@{=}[r]             &B,
}
$$
and
$$
\xyma{
\im(\gd)\ar@{>->}[r]^\iota\ar@{>->}[d] &B\ar@{->>}[r]\ar@{=}[d] &\coker(\gd)\ar[d]\\
B\ar@{=}[r]                            &B\ar[r]                 &0,
}
$$
one sees that there are short exact sequences of $\d$-bimodules of the form
$K'[1]\into\gd^\cdot\onto\iota^\cdot$ and $\iota^\cdot\into\co(B)\onto
K''[1]$, so that linear extensions by any $\gd$ can be described in terms of
those by bimodules of the form $K[1]$ and $\co(B)$. 

Before dealing with these, just let us make a note about lower cohomologies --
they can be expressed using derivations similarly to Hochschild cohomology. 

\dfn{der}
The group $\der(\d;\gd^\cdot)$ of \emph{derivations} of a linear form
$\d:M\to R$ with values in a $\d$-bimodule $\gd^\cdot=(B\ox_RM\to K\to B)$
consists of pairs of abelian group homomorphisms $(d:R\to B,\nabla:M\to K)$
satisfying 
\begin{align*}
&d\d=\gd\nabla,\\
&d(rs)=d(r)s+rd(s),\\
&\nabla(rm)=d(r)m+r\nabla(m),
\end{align*}
under pointwise addition. Its subgroup $\ider(\d;\gd^\cdot)$ consists of
\emph{inner} derivations $\ad(k)=(d_k,\nabla_k)$ for $k\in K$, defined by 
$$
d_k(r)=r\gd(k)-\gd(k)r,\ \ \nabla_k(m)=\d(m)k-\gd(k)\cdot m.
$$

\ 

Then by analogy with well known classical facts, 4.7 and 4.10 of \cite{jp}
readily gives 

\prp{lower}
{\sl
For a linear form $\d:M\to R$ and a $\d$-bimodule
$\gd^\cdot=(B\ox_RM\to K\to B)$, one has an exact sequence 
$$
0\to H^0(\d;\gd^\cdot)\to K\xto\ad{}\der(\d;\gd^\cdot)\to H^1(\d;\gd^\cdot)\to0.
$$
In other words, there are isomorphisms
$$
H^0(\t_\d;\gd^\cdot)\cong\setof{c\in K}{\forall m\in M\ (\d m)c=(\gd c)\cdot m}
$$ 
and}
$$
H^1(\t_\d;\gd^\cdot)\cong\der(\d;\gd^\cdot)/\ider(\d;\gd^\cdot).
$$ \qed

\ 

Now for $\co(B)$, one has

\prp{cone}
{\sl
For a linear form $\d:M\to R$ and an $R$-$R$-bimodule $B$, there is an
isomorphism 
$$
H^2(\t_\d;\co(B))\cong H^2(R;B),
$$
the latter being the MacLane cohomology group.}
\begin{proof}
Observe the diagram
$$
\xyma{
B\ar@{>->}[r]\ar@{=}[d] &N\ar@{->>}[r]\ar[d] &M\ar[d]^\d\\
B\ar@{>->}[r]           &S\ar@{->>}[r]       &R.
}
$$
It shows that the right hand square is pullback, so that the upper row is
completely determined by the lower one. Thus forgetting the upper row defines
an isomorphism, with the inverse which assigns to a singular extension of $R$
by $B$ the pullback as above. 
\end{proof}

Thus one arrives at a well studied situation here. As for the $K[1]$ case, we
have 

\prp{one}
{\sl
For a linear form $\d:M\to R$ and a left $R$-module $K$, there is an
isomorphism}
$$
H^2(\t_\d;K[1])\cong\ext^1_R(M,K).
$$
\begin{proof}
This is obvious from the diagram
$$
\xyma{
K\ar@{>->}[r]\ar[d] &{*}\ar@{->>}[r]^p\ar[d]^{\d p} &M\ar[d]^\d\\
0\ar[r]             &R\ar@{=}[r]                    &R.
}
$$
\end{proof}

But moreover the diagram
$$
\xyma{
\ker(\d)\ar@{>->}[r]\ar[d] &M\ar@{->>}[r]\ar[d]^\d &\im(\d)\ar@{>->}[d]\\
0\ar[r]                    &R\ar@{=}[r]            &R
}
$$
shows that $\t_\d$ is itself a linear extension of a theory corresponding to
the linear form of type $\gth a\into R$, where $\gth a$ is a left
ideal in $R$, by a natural system corresponding to an $(\gth a\into
R)$-bimodule of the form $K[1]$. And this is clearly the end: one obviously
has 

\prp{end}
{\sl
An abelian theory without constants cannot be represented nontrivially as a
linear extension of another theory if and only if it is of the type
$\t_{\gth a\into R}$, for the left linear form determined by a left ideal
$\gth a$ in a ring $R$ which does not have any nontrivial square zero
two-sided ideals.}
\begin{proof}
The only nontrivial remark to make here is that for any square zero two-sided
ideal $\gth b\into R$, one gets an extension 
$$
\xyma{
\gth k\ar@{>->}[r]\ar@{>->}[d] &\gth a\ar@{->>}[r]\ar@{>->}[d] &\gth a/\gth k\ar[d]\\
\gth b\ar@{>->}[r]             &R\ar@{->>}[r]               &R/\gth b
}
$$
for any left ideal $\gth k$ with
$\gth{ba}\subseteq\gth k\subseteq\gth b\cap\gth a$.
\end{proof}

Finally, consider an abelian theory $\t$ with constants. It has a largest
subtheory $\t_0$ without constants, obtained from $\t$ by removing all
morphisms $1\to X^n$ for $n>0$. The constants of $\t$ will then reappear in
$\t_0$ as \emph{pseudoconstants}, that is, those unary operations $p:X\to X$
satisfying the identity $p(a)=p(b)$. Conversely, if a theory without constants
is obtained nontrivially in such way, it must have some pseudoconstants. 

We know by \ref{thm:crab} that $\t_0=\t_\d$ for some linear form $\d:M\to R$.
Now pseudoconstants in $\t_\d$ correspond to elements $p$ of $M$ satisfying
the identity $\ph_a(p)=\ph_b(p)$ (with $a,b$ as variables). Using the identity
$(*)$ from \ref{dfn:aff} this gives $(1-\d p)_ab=a$, i.~e. $(\d p)_ab=b$ for
any $a$, $b$ in any affinity. Then taking $a=\brk{0,0}$, $b=\brk{0,1}$ in
$M\x R$ gives $\d p=1$. Now clearly there is a $p\in M$ with $\d p=1$ if and
only if $\d$ is surjective, in which case it is split by $\sigma(r)=rp$. Thus
in this case our linear form is isomorphic to 
(projection)$:\ker(\d)\os R\to R$. Let us fix one such $p$. We then may
declare $\setof{p+t_0}{t_0\in K}$ to be the set of pseudoconstants
corresponding to nullary operations, where $K$ is either empty or any
$R$-submodule of $\ker(\d)$. All choices will give equivalent categories of
models, the only difference being that for $K=\varnothing$ the empty set is
also allowed as a model. Each other model $A$ shall then have at least one
element $a$, and value of the unary operation $p$ on $A$ at $a$ will then be
$\ph_a(p)$, which does not depend on $a$ as we just saw. Denoting this element
by $0_A$ fixes a canonical $R$-module structure on each non-empty model.
Moreover each element of $M$ becomes uniquely written as $m=k+rp$ with
$r\in R$ and $k\in\ker(\d)$, so $\ph_a(m)=\ph_a(k)+r$, i.~e. $\ph_a$ is
completely determined by its restriction to $\ker(\d)$. Moreover by (*) of
\ref{dfn:aff} it is determined by $\ph_{0_A}$ alone. We see that, ignoring the
possible empty model, the category $\t_\d\m$ is equivalent to the coslice
$(\ker\d)/R\m$. We thus have proved: 

\prp{lint}
{\sl
An abelian theory has at least one pseudoconstant if and only if the category
of its models is equivalent to the category $K/(R\m)$ of left $R$-modules
under $K$, for some ring $R$ and an $R$-module $K$, with the possible
difference that the empty set is another model.
}\qed

\ 

Now observing \ref{exm}.\ref{mod} we conclude

\crl{lint}
{\sl
Any abelian theory with constants is isomorphic to $\t_{R;K}$ (defined in
\ref{exm}.\ref{mod}) for some ring $R$ and left $R$-module $K$.}\qed

\ 

Concerning linear extensions one observes that by \ref{exm}.\ref{mod}, any
theory with constants $\t_{R;K}$ is a trivial untwisted linear extension
of $\t_R$ by the bifunctor constructed there. Also observe that in any
linear extension $\t'\to\t$ of abelian theories one has constants if and
only if the other does. 

\ 

On the other hand, a description similar to \ref{prp:lint} is in fact possible
for categories of models of abelian theories without constants too. For any
left linear form $\d:M\to R$, denote (temporarily) by $\d\aff'$ the following
category: objects are $R$-module homomorphisms $f:M\to N$; a morphism from
$f':M\to N'$ to $f:M\to N$ is a pair $(g,n)$, where $g:N'\to N$ is an
$R$-module homomorphism and $n\in N$ is an element such that
$f(x)-gf'(x)=\d(x)n$ holds for all $x\in M$. Composition is given by
$(g,n)(g',n')=(gg',n+g(n'))$, and identities have form $(\id,0)$.
Equivalently, one might define objects as commutative triangles 
$$
\xyma{
&M\ar[dr]\ar[dd]_\d\\
\ \ \ \ \ \ \ \ \ &&R\os N\ar[dl]^{\textrm{projection}}\\
&R
}
$$
and morphisms as commutative diagrams
$$
\xyma{
&M\ar[dl]\ar[dr]\\
R\os N'\ar[dr]\ar[rr]&&R\os N\ar[dl]\\
&R
}
$$
in $R\m$. One then has

\prp{varaff}
{\sl
The category $\t_\d\m$ is equivalent to $\d\aff'$ with an extra initial object
added.}
\begin{proof}
Define a functor $\d\aff'\to\d\aff$ as follows: for $f:M\to N$, define an
affinity structure on $N$ by $a+_bc=a-b+c$, $r_ab=(1-r)a+rb$, and
$\ph_a(x)=f(x)+(1-\d x)a$. And to a morphism $(g,n)$ assign the homomorphism
of affinities $N'\to N$ given by $n'\mapsto n+g(n')$. It is straightforward to
check that this defines a full and faithful functor. Moreover any nonempty
affinity is isomorphic to one in the image of this functor -- just choose an
element and use it as zero to define a module structure and a homomorphism
from $M$ according to the affinity identities.
\end{proof}

This allows to give an example, which looks pleasantly familiar:

\subsection*{Example} Fix a field $k$, and let the \emph{category of cycles}
be defined as follows. Objects are pairs $((V,d),c)$, where $(V,d)$ is a
differential $k$-vector space and $c\in V$ is a cycle, i.~e. $dc=0$. A
morphism from $((V,d),c)$ to $((V',d'),c')$ is a pair $(\ph,x)$, where
$\ph:V\to V'$ is a $k$-linear differential map and $x\in V'$ an element with
$c'-\ph(c)=dx$. With the evident identities and composition this forms a
category which is clearly of the form $\t_\d\m$, for the linear form
$\d:\eps k[\eps]\into k[\eps]$, where $\eps$ is an indeterminate element with
$\eps^2=0$. 

Now obviously this example admits a linear extension structure over $\t_k$,
since $\eps k[\eps]$ is a square zero ideal. But of course \ref{prp:end}
provides lots of similar (less cute) examples without this property.

\bibliographystyle{amsalpha}

\end{document}